\documentclass[letterpaper,oneside]{amsart}
\pdfoutput=1 
\usepackage[T1]{fontenc}
\usepackage{lmodern} 
\usepackage[defaultsans]{cantarell} 
\usepackage{color,soul}
\usepackage{amsrefs}
\usepackage[pagewise]{lineno}
\usepackage[dvipsnames]{xcolor}
\usepackage{todonotes}
\setuptodonotes{inline}
\usepackage[utf8]{inputenc} 
\usepackage{comment}
\usepackage{mathtools} 
\usepackage{colonequals}
\usepackage{amsmath}
\usepackage{amssymb}
\usepackage{tikz}
\usepackage{mathrsfs}
\usepackage{enumitem}
\usepackage{mathtools}
\usepackage{amsaddr}

\usepackage{mleftright}
\usepackage{graphicx}
\usepackage{hyperref}
\usepackage{cancel}
\usepackage[english]{babel} 

\usepackage{tikz}
\usepackage{thmtools}
\usepackage{amssymb}
\usepackage[all]{xy} 
\usepackage{xparse}
\usepackage[notocbasic,refpage]{nomencl}
\makenomenclature

\newcommand{\Id}[1]{\operatorname{Id}_{#1}} 
\newcommand{\R}{\mathbb{R}}
\newcommand{\Z}{\mathbb{Z}}

\newcommand{\catname}[1]{\mathbf{#1}}

\newcommand{\Absurd}{(\rightarrow\leftarrow)}

\newcommand{\GR}[1]{\left| #1 \right|} 
\newcommand{\CC}[1]{\mathbf{K}\!\left(#1 \right)} 
\newcommand{\GRCC}[1]{\GR{\CC{#1}}} 
\newcommand{\uGRCC}[1]{\operatorname{U}\left|#1\right|} 

\newcommand{\GFil}[2]{\left[#1 \right]_{#2}} 
\newcommand{\iGR}[1]{\mathbf{R}\!\left( #1 \right)} 
\NewDocumentCommand{\DMod}{m o}{ 
	\Theta_{#1}\IfValueT{#2}{\!\left({#2} \right)} 
}
\newcommand{\Fld}[3]{\left( #1  \right)_{#2}^{#3}} 
\newcommand{\FldF}[2]{#1^{#2}} 
\newcommand{\ConM}[2]{\left. #1\right._{#2}^{\text{\textbf{conv}}}} 
\newcommand{\Fil}[1]{\mathtt{Fil}\!\left( #1\right)} 
\newcommand{\FilU}[1]{\mathtt{UFil}\!\left( #1\right)} 
\newcommand{\NFil}[1]{\mathcal{U}_{#1}} 
\newcommand{\PowS}[1]{2^{#1}} 
\newcommand{\ResF}[2]{#1\mid_{#2}} 

\newcommand{\Ball}[3]{\operatorname{B}_{#3}(#1,#2)} 

\title[Graphs and Their Clique Complexes]{Graphs and Their Clique Complexes Are Pseudotopological Weakly Homotopy Equivalent}

\author{Jonathan Treviño-Marroquín$^{1}$}
\address{$^1$Universidad Autónoma de Nuevo León, Nuevo León, México}
\email{jonathan.trevinomrqn@uanl.edu.mx}

\thanks{This work was supported by the CONACYT, now SECIHTI, postgraduate studies scholarship number 839062, and the results in this article form a part of the author's PhD thesis project supervised by Antonio Rieser at the Centro de Investigación en Matemáticas, Guanajuato, México.}

\theoremstyle{plain}
\newtheorem{teorema}{Theorem}[section]
\newtheorem{proposicion}[teorema]{Proposition}

\newtheorem{lema}[teorema]{Lemma}

\theoremstyle{definition}
\newtheorem{definicion}[teorema]{Definition}

\newtheorem{observacion}[teorema]{Remark}
\newtheorem{ejemplo}[teorema]{Example}

\declaretheoremstyle[
qed=\qedsymbol
]{mystyle}

\begin{document}

\begin{abstract}
We establish a connection between the graphs and their geometric realizations of their clique complexes, proving that graphs with their canonical pseudotopological structure, are weakly homotopy equivalent to the geometric realization of their clique complex.
\end{abstract}

\maketitle





\section*{Introduction}

\pagenumbering{arabic}\label{cap.introduccion}

	Recently, there has been increased activity and interest in the homotopical study of reflexive graphs from several different perspectives.
	For example, in \cite{Babson_etal_2006} the authors develop the $A$-theory in which cubical-type homotopy is defined in graphs using the inductive product and the discrete interval to generate their homotopies.
	In \cites{Dochtermann_2009, Chih_Scull_2021} $\times$-homotopy is introduced to study not-necessarily reflexive graphs.
	There are also homotopy theories for digraphs in \cite{Grigoryan_etal_2014} and the non-topological homotopy for metric spaces considered in \cite{Plaut_Wilkins_2013}.
	Restricting to reflexive graphs on cubical lattices, where the edges are determined by their position in the lattice, yields digital homotopy \cites{Lupton_etal_2022a,Lupton_Scoville_2022}.
	Finally, in \cites{Rieser_2021,Rieser_arXiv_2022} Antonio Rieser proposed using \v{C}ech closure spaces and pseudotopological spaces as categories in which to develop the homotopy theories of topological spaces, graphs, digraphs simultaneously, and then in point clouds (the latter if which is used in topological data analysis \cites{Carlsson_2009, Carlsson_et_al_2012}), thereby extending the earlier work of Demaria \cite{Demaria_1987} on graphs. 

	In this article, we continue developing homotopy theory in pseudotopological spaces, as Rieser begun in \cites{Rieser_2021,Rieser_arXiv_2022}.
	We prove that every graph $X$ is weakly homotopy equivalent to its clique complex $\GRCC{X}$.
	To achieve this, we define an intermediate pretopological space $\iGR{X}$ (\autoref{def:GraphRealizations}), which is a simplicial complex with a modified closure structure.
	This structure ensures such that (1) there is a continuous inclusion from the vertices of $X$ into $\iGR{X}$, (2) $\iGR{X}$ is strongly homotopy equivalent to the $X$, and (3) $\iGR{X}$ is weakly homotopy equivalent to the geometric realization of the clique complex of $X$, denoted by $\GRCC{X}$. 
	The construction of the weak homotopy equivalence further requires us to construct three transformations of functions $f:Y \to \iGR{X}$, which are the key to the proof.
\begin{itemize}
	\item \textbf{The discrete modification} (\autoref{def:DiscModification}), which sends a map $f$ with codomain $\iGR{X}$ to a map whose image lies entirely within $X \subset \iGR{X}$.

	\item \textbf{The flood} (\autoref{def:Flood}), which expands the preimage of a vertex $v \in X \subset \iGR{X}$ to create a neighborhood $U$ of $f^{-1}(v) \subset Y$ that is small enough so that the image of every element in $U$ is contained in the closure of the image of any of its points.

	\item \textbf{The convex transformation} (\autoref{def:ConvexTransformation}), which converts maps $f:Y \to \iGR{X}$ whose image is in $X \subset \iGR{X}$ into a simplicial map from a triangulation of $Y$ to $\iGR{X}$.
\end{itemize}
	We study the properties of each of these transformations in turn, and then we use in combination to prove our main result.
	The critical property of each transformation is that, when applied to a continuous map from a triangularizable compact metric space, they produce a continuous map.

	Using the flood transformation after the discrete modification employs the same fundamental idea as that studied in digital homotopy to determine the second homotopy group for the digital sphere \cite{Lupton2023DigitalHom}.
	Unlike that paper, which uses the lattice $\Z^2$ as the domain and discrete distances, this paper uses a triangularizable compact metric space as its domain.
	More specifically, it uses $S^n$ and $S^n\times I$. 

	When combined with well-known results about the Vietoris-Rips complex \cite{Latschev_2001}, the results in this paper show that, for a ``big enough'' cloud of points sampled from a manifold, the pseudotopological homotopy groups (defined in \cite{Rieser_arXiv_2022}) from certain closure spaces built on those points are isomorphic to the classical homotopy groups of the manifold being sampled.
	The Vietoris-Rips complex is difficult to understand and only a few complete calculations are known for metric spaces at any scale, for example \cites{Adams2023lower, Adams2022VietorisRipsHypercube, Adams2020VietorisRipsMetricGluings} and is used to bound Gromov-Hausdorff distance \cites{Adams2022gromovhausdorff, Memoli2023GHUltraMetricSpaces, Memoli2023GHSpheres}.
	Thus, with this paper and growing wave of studies in these categories mentioned above, we may hope to find new ways to approach the Vietoris-Rips complex.

	The isomorphism between the pseudotopological homotopy groups of finite graphs and those of their clique complexes was also independently established by Mili\'cevi\'c and Scoville in \cite{WebJMMNicks}, later on \cite{MilicevicScoville2026}.
	Their proof differs from this one, particularly as it does not pass through the space $\iGR{X}$. Additionally, their proof applies to finite digraphs, whereas the technique we study here is also valid for arbitrary, even uncountable infinite, reflexive graphs.
	
	We have a discussion section at the end of the document where we compare the results to some of those in \cite{Adamaszek_Adams_2017} and \cite{MilicevicScoville2026}.
	Additionally, we added a nomenclature section to make the notation used easier to identify.

\section{Graphs and its Geometric Realizations}
\label{cap.WHE}	

This paper discus objects and morphisms in the category of pseudotopological spaces, particularly, pretopological ones.
The latter category is related to the category of \v{C}ech closure spaces.
Proposition 2.3.1.8 \cite{Preuss_2002} provides the full argument proving these two categories are equivalent. We recall the following definitions:
	
\begin{definicion}
	\label{def:FilterGenerated}
	For every set $X$ and a collection $\mathcal{A}$ of subsets of $X$ with the finite intersection property, \emph{filter generated by $\mathcal{A}$} is defined as follows:
	\[\GFil{\mathcal{A}}{X} = \{Y\subset X \mid Y\supset A \text{ for a }A\in\mathcal{A}\}.\]\nomenclature{$\GFil{\mathcal{A}}{X}$}{The filter generated by $\mathcal{A}$ on $X$}
	Whenever $\mathcal{A} = \{ A \}$, we write $\GFil{A}{X}$ and if $A=\{x_0\}$, we use the notation $\GFil{x_0}{X}$.
	
	The collection of all the filters on $X$ is denoted by $\Fil{X}$\nomenclature{$\Fil{X}$}{The collection of filters on $X$} and the collection of all the ultrafilters on $X$ is denoted by $\FilU{X}$.\nomenclature{$\FilU{X}$}{The collection of ultrafilters on $X$}
\end{definicion}

\begin{definicion}[\cite{Beattie_Butzmann_2002}, 1.1.1 and 1.3.23] 
	\label{def:ClosureSpace}
	A \emph{pseudotopological space} is a pair $(X,\lambda)$, where $X$ is a set and $\lambda:X\rightarrow \PowS{\Fil{X}}$ a function that satisfies following three conditions:
	\begin{enumerate}
	\item $\GFil{x}{X}$ belongs to $\lambda(x)$.
	
	\item For all filters $\mathcal{F}$ and $\mathcal{G}$ in $\lambda(x)$, the set intersection $\mathcal{F}\cap \mathcal{G}$ belongs to $\lambda(x)$.
	
	\item If $\mathcal{F}$ belongs to $\lambda(x)$, so does every filter $\mathcal{G}$ finer than $\mathcal{F}$.
	
	\item $\mathcal{F}$ belongs to $\lambda(x)$ whenever every ultrafilter $\mathcal{G}$ finer than $\mathcal{F}$ belongs to $\mathcal{\lambda}(x)$.
	\end{enumerate}

	Additionally, a pseudotopological space is called a pretopological space if
	
	\begin{enumerate}\setcounter{enumi}{4}
		\item The intersection of all the filters in $\lambda(x)$ belongs to $\lambda(x)$.
	\end{enumerate}
\end{definicion}

	The intersection of all the filters in $\lambda(x)$ is called the \emph{neighborhood filter at $x$ in $X$} and is denoted $\NFil{x}$.\nomenclature{$\NFil{x}$}{The neighborhood filter at the point $x$}
	When $\mathcal{F}$ belongs to $\lambda(x)$, it is said that $\mathcal{F}$ converges to $x$ and it is denoted $\mathcal{F}\xrightarrow{\lambda} x$ or, simply, $\mathcal{F}\rightarrow x$.

\begin{definicion}\label{def:ContPsTop}
	Let $(X,\lambda)$ and $(Y,\mu)$ be pseudotopological spaces.
	A map $f:X\rightarrow Y$ is called continuous from $(X,\lambda)$ to $(Y,\mu)$, denoted as $f:(X,\lambda)\rightarrow (Y,\mu)$, if $f(\mathcal{F})\in \mu(f(x))$ for every $\mathcal{F} \in \lambda(x)$ and for every $x\in X$.
	
	It is said that $f:X\to Y$ is continuous at $x\in X$ if $f(\mathcal{F})\in \mu(f(x))$ whenever $\mathcal{F}\in \lambda(x)$. It can be proven that $f:X\to Y$ is continuous if and only if $f$ is continuous at $x$ for every $x\in X$.
\end{definicion}

	When the context is clear, we omit the references to the spaces and we say that $f$ is continuous. Note that if $(X,\lambda)$ and $(Y,\mu)$ are pretopological spaces, then $f:X\to Y$ is continuous if and only if $f(\NFil{x}) \supset \NFil{f(x)}$ for every $x\in X$.

	We denote the category whose objects are pseudotopological spaces and whose morphisms are continuous maps by $\catname{PsTop}$. The category whose objects are pretopological spaces and whose morphisms are continuous maps $\catname{PrTop}$.
	
\begin{observacion}
	Every pseudotopological space $(X,\lambda)$ induces a \v{C}ech closure space, $(X,c_\lambda)$, which is defined using the unique closure operator induced for the neighborhood filters $\NFil{x}$ for every $x\in X$. When $(X,\lambda)$ is a pretopological space, then this relation is an equivalence of categories.
	
	More specifically:
	\begin{itemize}
		\item If $(X,c)$ is a \v{C}ech closure space, then consider $\NFil{x}$ as the collection of all the sets $U\subset X$ such that $x \in X-c(X-U)$.
		
		\item If $(X,\lambda)$ is a pretopological space, then  $c:\PowS{X}\to \PowS{X}$ is defined as \[c(A) = \{ x\in X~\mid~\forall U\in\NFil{x},~U\cap A\neq \varnothing \}\}\]
		for every $A\subset X$. 
	\end{itemize}
	To see more about this, read \cite{Rieser_2021}.
\end{observacion}
	
\begin{ejemplo}
	Every topological spaces can be seen as a pretopological space with the neighborhood filters.
\end{ejemplo}

	The following is an example of a continuous map from a topological space to a non-topological space.

\begin{ejemplo}
	Let $S^1 = \{ z\in \mathbb{C} \mid~|z|=1  \} $ with the euclidean topology, and let $C_4$ be a 4-cyclic graph (or the digital circle), i.e., the set $\{0,1,2,3\}$ with the closure $cl(y)=\{y-1,y,y+1\}$ module 4 defined for every $y$ in the set.
	We have that the map $f:S^1\rightarrow C_4$ defined by
\begin{align*}
f(x) = \left\lbrace\begin{array}{ll}
0 & \text{ if } x \in [0,1/4);\\
1 & \text{ if } x \in [1/4,2/4);\\
2 & \text{ if } x \in [2/4,3/4);\\
3 & \text{ if } x \in [3/4,1).
\end{array}\right.
\end{align*}
	is continuous.
\end{ejemplo}

In \cite{Rieser_arXiv_2022}, Rieser defines and studies the pseudotopological homotopy and the weak homotopy equivalence for $\catname{PsTop}$.

\begin{definicion}
	\label{def:Homotopy}
	Let $f,g: X\rightarrow Y$ be continuous maps between pseudotopological spaces. The maps $f$ and $g$ are said to be \emph{homotopic}, denoted by $f\simeq g$, if there exists a continuous map $H:X\times I\rightarrow Y$ such that $H(x,0)=f(x)$ and $H(x,1)=g(x)$ for every $x\in X$.
\end{definicion}

\begin{definicion}
	\label{def:SHE}
	The pseudotopological spaces $X$ and $Y$ are called \emph{strongly homotopy equivalent} if there exist continuous maps $f: X\rightarrow Y$ and $g: Y\rightarrow X$ such that $gf\simeq \Id{X}$ and $fg\simeq \Id{Y}$.
\end{definicion}

\begin{definicion}
	\label{def:HomotopyClasses}
	For every natural number $n$, we define the \emph{homotopy classes} of a pseudotopological space $X$ as the set of equivalence classes of continuous maps $(S^n,*)\rightarrow (X,*)$.
	We denote this set by $\pi_n(X,*)$.
\end{definicion}
	
	Using the standard $\star$ operation on these maps and standard arguments, we can prove that $\pi_n(X,*)$ is a group for all $n\geq 1$ and an abelian group for all $n>1$.
	
\begin{definicion}
	\label{def:WeakHomEq}
	Let $f:X\to Y$ be a continuous map. The map $f$ is \emph{weak homotopy equivalence}, if $\pi_n(f): \pi_n(X,*)\to \pi_n(Y,*)$ is an isomorphism for every $n\geq 0$. The pseudotopological spaces $X$ and $Y$ are \emph{weakly homotopy equivalent} if there exist a weak homotopy equivalences $f:X\to Y$ or $g:Y\to X$.
	
	We say that $X$ and $Y$ have the same weak homotopy type if a zig-zag of weak homotopy equivalences exists from $X$ to $Y$
\begin{align*}
X \rightarrow X_1 \leftarrow X_2 \rightarrow \ldots \leftarrow X_{n-1} \rightarrow X_{n} \leftarrow Y.
\end{align*}
\end{definicion}

\begin{definicion}\label{def:CanStrGraph}
	Let $(X, E)$ be a simple, reflexive graph, meaning the graph has no multiple edges and $(x,x)$ is an edge for every $x\in X$.
	
	The \emph{canonical closure operator for a graph $(X,E)$}, denoted $c_X$, is given by
\begin{align*}
c_X(x) \coloneqq & \{ y\in X \mid (x,y)\in E \} \text{ for every }x\in X, \\
c_X(A) \coloneqq & \bigcup\limits_{x\in A} c_X(A)  \text{ for every }A\subset X.
\end{align*}
\end{definicion}

	Note that this closure operation induces the neighborhood filter at $x$ on $X$ as $\NFil{x} = \GFil{c_X(x)}{X}$.
	When the context is clear, we omit writing the set of edges in the ordered pair. In other words, we write $X$ instead of $(X, E)$ when we refer to a graph.
	
\begin{observacion}
	For any set, we can prove the existence of a well-order through the axiom of choice.
	The set of vertices in each of our graphs is a well-ordered.
	This order is denoted by $\leq$.
	The symbols $<$ will mean that $v\leq w$ and $v\neq w$.
	
	Our results holds even when considering different well-orders in the vertices.
	
	In the simplices, the well-order of the elements defines the order of the elements, not their appearance order. For example, when writing an arbitrary simplex $\{v_0,v_1,\ldots,v_n\}$, we understand that $v_0<v_1<\ldots<v_n$.
\end{observacion}

	The definition below introduces two geometric realizations of a graph.
	The first is the classical geometric realization induced by (the clique complex of) the graph.
	The second is a hybrid of the graph and the clique complex.
	The latter has a strong homotopy equivalence with the canonical pretopological structure of the graph and it is finer than the classical geometric realization.

\begin{definicion}
\label{def:GraphRealizations}
	Let $X$ be a graph.
\begin{enumerate}
	\item The (topological) geometric realization of the clique complex of $X$ is denoted by $\GRCC{X}$.
	We write the underlying set of $\GRCC{X}$ as $\uGRCC{X}$ and denote the closure operator by $c_{\GRCC{X}}$ and the pseudotopological spaces as $\GRCC{X} = (\uGRCC{X},\lambda_{\GRCC{X}})$.
	\nomenclature{$\CC{X}$}{The clique complex of the graph $X$}
	\nomenclature{$\GRCC{X}$}{The geometric realization of the clique complex of the graph $X$}
	\nomenclature{$\uGRCC{X}$}{The underlying structure of $\GRCC{X}$}
	
	\item The \emph{intermediate geometric realization} $\iGR{X}$ is defined over the same underlying set of $|\GRCC{X}|$, $\uGRCC{X}$, with the pretopological structure defined with the neighborhood filter at $v\in X$ as
	\begin{align*}
		\NFil{v}\coloneqq \GFil{\bigcup \{|\sigma|~\mid x\in \sigma,~\sigma\in\CC{X}\}}{\uGRCC{X}},
	\end{align*}
	\nomenclature{$\iGR{X}$}{The intermediate geometric realization of the graph $X$}
	and the neighborhood filter at $v$ on $\iGR{x}$ is equal to the one on $\GRCC{X}$ if $v\in \uGRCC{X}-X$.
	We denote the closure operator as $c_{\iGR{X}}$ and the pseudotopological spaces as $\iGR{X} = (\uGRCC{X}, \lambda_{\iGR{X}} )$.
\end{enumerate}
\end{definicion}

	The space $\iGR{X}$ and $\GRCC{X}$ have the same underlying set.
	By construction, there are fewer sets in the neighborhood filter at a vertex on $\iGR{X}$ than on $\GRCC{X}$, and they have the same neighborhood filter at $x$ when $x$ is not a vertex.
	Thus, $\iGR{X}$ is finer than $\GRCC{X}$, and the identity from $\GRCC{X}$ to $\iGR{X}$ is continuous.
	Similarly, $X$ is a (pretopological) space of $\iGR{X}$.

\begin{definicion}[The Barycentric Cover of a Simplex]
\label{def:BaryCoverSimplex}
Let $\Delta^n\coloneqq\{v_0,\ldots,v_n\}$ be an $n$-simplex. We denote the barycenter of $\{v_{i_0},\dots,v_{i_k}\}\subset \Delta^n$ as
\begin{align*}
b_{\{v_{i_0},\dots,v_{i_k}\}}.
\end{align*}
We define the sets
\begin{align*}
A_i \coloneqq \bigcup_{\alpha_j\in [n]} \{| \{ v_i,b_{\{v_i,v_{\alpha_1}\}}, \ldots, b_{\{v_i,v_{\alpha_1},\dots,v_{\alpha_n}\}}  \}|~\mid \alpha_j\in [n] \}.
\end{align*}
for every $i\in [n] = \{0,\ldots, n\}$.\nomenclature{$[n]$}{The set of integers $\{0,1,2,\ldots,n\}$}
The collection $\{A_0,\ldots,A_n\}$ is said to be \emph{the barycentric cover of $|\Delta^n|$}.
\end{definicion}

\begin{proposicion}
\label{prop:BaryCoverSimplex}
Let $\Delta^n\coloneqq \{v_0,\ldots,v_n\}$ be an $n$-simplex and $\{A_0,\ldots,A_n\}$ the barycentric cover of $\Delta^n$. Then $|\Delta^n|=\cup_{i=0}^n |A_i|$.
\end{proposicion}

\begin{proof}
Let $x\in |\Delta^n|$, then for every $i\in [n]$ there exists $t_i\geq 0$ such that $\sum_{i=0}^n t_i = 1$ and
\begin{align*}
x=\sum_{i=0}^n t_iv_i.
\end{align*}
Since we have a finite number of elements, there exists an index in $[n]$, say $\beta_0$, such that $t_{\beta_0}=\min_{i\in [n]}\{t_i\}$.
Thus
\begin{align*}
x= \sum_{i\in [n]- \{\beta_0\}} (t_i-t_{\beta_0}) x_i + \frac{t_{\beta_0}}{n+1}\left(\frac{1}{n+1} \sum_{i=0}^n x_i \right).
\end{align*}
Using the same argument, we can find an index $\beta_1$ such that $t_{\beta_1}=\min_{i\in [n]\backslash\beta_0}\{t_i\}$ and
\begin{align*}
x= & \sum_{i\in [n]\backslash \{\beta_0,\beta_1\}} (t_i-t_{\beta_1}) x_i + (t_{\beta_1}-t_{\beta_0})n\left(\frac{1}{n} \sum_{i\in [n]\backslash\{\beta_0\}}x_i \right)\\
&+ t_{\beta_0}(n+1)\left(\frac{1}{n+1} \sum_{i=0}^n x_i \right).
\end{align*}
By following this procedure inductively until all vertices are covered, we obtain that
\begin{align*}
x = \sum_{k=0}^n \left( (t_{\beta_{n-k}}-t_{\beta_{n-k-1}})(k+1)\left( \frac{1}{k+1} \sum_{i\in [n]\backslash \{\beta_j\mid j\in [k] \}}x_i \right)  \right),
\end{align*}
	considering $t_{\beta_{-1}}=0$ to simplify the expression.
	Note that, by definition of barycenter,
\begin{align*}
b_{\{v_{\beta_{j}}\mid j\in [k]\}} = \frac{1}{k+1} \sum_{i\in [n]\backslash \{\beta_j\mid j\in [k] \}}x_i.
\end{align*}
Concluding that $x\in A_{\beta_{n}}$.
\end{proof}

Take a $n$-simplex $\Delta^n\coloneqq \{v_0,\ldots,v_n\}$ and a $m$-simplex $\Delta^m\coloneqq \{w_0,\ldots,w_m\}$. If $\Delta^n$ and $\Delta^m$ meets, then its intersection is a $k$-simplex, say $\Delta^k\coloneqq \{v_0,\dots,v_k\}$ with $v_i=w_i$ for every $i\in [k]$.
Call $\{A_0,\ldots,A_n\}$, $\{B_0,\ldots,B_m\}$ and $\{C_0,\ldots, C_k\}$ the barycentric cover of $\Delta^n$, $\Delta^m$ and $\Delta^k$, respectively.
Note that \[C_i=A_i\cap |\Delta^k| = B_i\cap |\Delta^k|\] for every $i\in [k]$ given that the barycentric cover of $C_i$ only depends on the barycenters in $\Delta^k$. Given this argument, we can inductively extend the following definition (\autoref{def:DiscModification}) to all $\iGR{X}$ applying it to every simplex of $\CC{X}$.

\begin{definicion}[Discrete Modification]
\label{def:DiscModification}
	Let $X$ be a graph. \emph{The discrete modification} $\DMod{X}:\iGR{X}\rightarrow X\subset \iGR{X}$ is defined as follows:
	\nomenclature{$\DMod{X}$}{The discrete modification of $\iGR{X}$}
	Let $\Delta^n\coloneqq \{v_0,\ldots,v_n\}\in \CC{X}$ be an $n$-simplex. The map $\ResF{\DMod{X}}{\Delta^n}: |\Delta^n|\to |\Delta^n|$ is the relation
	\[\ResF{\DMod{X}}{\Delta^n}(x)\coloneqq  \left\lbrace\begin{array}{ll}
v_0 & x\in A_0,\\
v_1 & x\in A_1-A_0,\\
& \vdots\\
v_n & x\in A_n-\cup_{i=0}^{n-1}A_i
\end{array}\right. .\]

The map $\DMod{X}:\uGRCC{X}\to X$ is defined such that $\DMod{X}(x) = \ResF{\DMod{X}}{\Delta^n}(x)$ for some $\Delta^n$ to which $x$ belongs. 
\end{definicion}

\begin{figure}[h]
\begin{tikzpicture}
\draw (3,0) -- (5,0) -- (4,1.732)--cycle;
\filldraw[Green] (3,0) circle (0.05);
\filldraw[Blue] (5,0) circle (0.05);
\filldraw[Red] (4,1.732) circle (0.05);

\draw[->] (2.25,1)--(2.75,1);

\draw (0,0) -- (2,0) -- (1,1.732)--cycle;
\filldraw[Blue] (2,0) -- (1,0) -- (1,.577)--(3/2,1.732/2)--cycle;
\filldraw[Green] (0,0) -- (1,0) -- (1,.577)--(1/2,1.732/2)--cycle;
\filldraw[Red] (1,1.732) -- (1/2,1.732/2) -- (1,.577)--(3/2,1.732/2)--cycle;
\end{tikzpicture}
\caption{Illustration of the discrete modification for the geometric realization of a 2-simplex.}
\label{fig:DiscreteModification}
\end{figure}
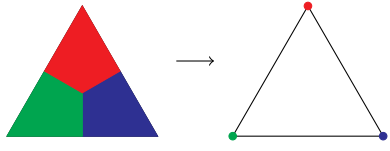

\begin{lema}
\label{lem:ThetaXContinuous}
	Let $X$ be a graph, the discrete modification is continuous.
\end{lema}

\begin{proof}
	Consider the following two case of $x\in \uGRCC{X}$. In the first case, let $x\in X$. Remember that the neighborhood filter at $x$ in $\iGR{X}$ is \[\NFil{x} = \GFil{\bigcup \{|\sigma|~\mid x\in \sigma,~\sigma\in\CC{X}\}}{\uGRCC{X}}.\]
	By the definition of the discrete modification, $\DMod{X} (|\sigma|) = \{v\mid v\in \sigma \}\subset |\sigma|$. Then \[\DMod{X}(\cup_{x\in \sigma} |\sigma|)=\cup_{x\in \sigma} \{ v \mid v\in\sigma \}\subset \cup_{x\in \sigma} |\sigma|.\] Thus $\DMod{X}(\GFil{\cup_{x\in \sigma} |\sigma|}{\uGRCC{X}})\rightarrow x = \DMod{X}(x)$.

	In the second case, let $x\notin X$, by the definition of $\uGRCC{X}$, there exists a unique simplex $\sigma$ containing $x$ in its relative interior.
	Then $\DMod{X}(x)\in \{v \mid v\in \sigma\}$.
	A filter $\mathcal{F}$ converges to $x$ in $\GRCC{X}$ (and then in $\iGR{X}$) if it converges in the restriction of every simplex containing $x$.
	This implies that $\cup_{\{v\mid v\in\sigma \}\subset\sigma'}|\sigma'| \in \NFil{x}$.
	Then \[\DMod{X}(\NFil{x})\supset \DMod{X} \left(\GFil{\cup_{\{v\mid v\in\sigma \}\subset\sigma'}|\sigma'|}{\uGRCC{X}}\right) \supset \GFil{\cup_{\DMod{X}(x)\in\sigma}|\sigma|}{\uGRCC{X}}.\]Thus $\DMod{X}(\NFil{X})\rightarrow \DMod{X}(x)$.
\end{proof}

\begin{observacion}
By construction, we have that $Im(\DMod{X})\subset X$, $\DMod{X} \circ \DMod{X} = \DMod{X}$ and $\DMod{X}(v)=v$ for every $v\in X\subset \uGRCC{X}$.
\end{observacion}

\begin{proposicion}
\label{prop:SHEBetweenGraphandRGraph}
Let $X$ be a graph. $X$ and $\iGR{X}$ are strongly homotopy equivalent.
\end{proposicion}

\begin{proof}
Let $g: X\rightarrow \iGR{X}$ be the map that sends the points in $X$ to their corresponding vertices of the simplices in $\uGRCC{X}$.
By definition, $\DMod{X}\circ g=\Id{X}$, so we only need to find a homotopy from $g\circ \DMod{X}$ to $\Id{\iGR{X}}$. 

First we restrict to a simplex $\Delta\in \CC{X}$. Let $\ResF{H}{\Delta}:|\Delta|\times I\rightarrow |\Delta|\subset R(X)$ be defined by
\begin{align*}
\ResF{H}{\Delta}(x,t) \coloneqq \left\lbrace \begin{array}{ll}
\Id{|\Delta|}(x) & \text{ if }t\in [0,1),\\
\DMod{X}(x) & \text{ if }t=1.
\end{array}\right.
\end{align*}
We claim that $\ResF{H}{\Delta}$ is continuous. For every point $(x,1)$, its image contained in $\Delta$, then $\ResF{H}{\Delta}$ is continuous at that point by the definition of the neighborhood filter at a vertex in $\iGR{X}$. For every point $(x,t)$ such that $t<1$, there exists $r>0$ such that the ball centered in $(x,t)$ with radius $r$, $\Ball{(x,t)}{r}{|\Delta|\times I}$, does not intersect the set $\Delta\times\{1\}$. This conclude that $H(\Ball{(x,t)}{r}{|\Delta|\times I}) = \Ball{x}{r}{|\Delta|}$ and $\ResF{H}{\Delta}$ is continuous at that point.

Define $H$ such that $H(x,t) = \ResF{H}{\Delta}(x,t)$ is $x\in |\Delta|$. Denote the neighborhood filter of $t$ in $[0,1]$ by $\mathcal{I}(t)$. For $t<1$ and $x\in \uGRCC{X}$, there exists an $\varepsilon>0$ such that for every $U\in\NFil{x}$ satisfies that $U\times (t-\varepsilon, t+\varepsilon)$ does not intersect $U\times\{1\}$.
Then, $H(\NFil{x}\times (t-\varepsilon,t+\varepsilon))=\NFil{x}$. Therefore, $H(\NFil{x}\times\mathcal{I}(t))\rightarrow x$.

On the other hand, when $t=1$,
\begin{align*}
	H(\NFil{x}\times\mathcal{I}(1)) = & \GFil{\{H(A\times [1-1/n,1)\mid A\in\NFil{x},n\in\mathbb{N}\}}{\uGRCC{X}}\\
	= & \GFil{\{\DMod{X}(A)\cup A\mid A\in\NFil{x}\}}{\uGRCC{X}}
\end{align*} \[\] by the definition of $H$. According to \autoref{lem:ThetaXContinuous}, there exists an element $B\in \NFil{x}$ such that $\DMod{X}(B)\subset \cup_{\DMod{X}(x)\in\sigma}|\sigma|$ and $B\subset \cup_{\DMod{X}(x)\in\sigma}|\sigma|$ for every $B\subset A$. Thus, $H(\NFil{x}\times\mathcal{I}(1))\subset [ \cup_{\DMod{X}(x)\in\sigma}|\sigma|]$ concluding that $H(\NFil{x}\times\mathcal{I}(1))\rightarrow \DMod{X}(x)$.
\end{proof}

	The following results will be used for maps on an $n$-sphere.
	However, they will be proven in metric spaces with certain properties that make them suitable for a more general context.
	From this point on, we set $Y$ as a metric space $(Y,d)$.
	In addition, for any metric space $Z$, we denote the open ball of radius r centered at $z$ in $Z$ by the set $\Ball{z}{r}{Z}\coloneqq \{ z'\in Z \mid d(z,z') <r \}$.
	
	Before the following lemma, we recall the next definition:
	
\begin{definicion}\label{def:LocFinGraph}
	Let $(X, E)$ be a graph. It is said that $(X,E)$ is \emph{locally finite} if for every $x\in X$ there exists a finite number of vertices $(x,y)\in E$.
\end{definicion}

\begin{lema}
\label{lem:DiscretizationOfAMap}
Let $X$ be a graph, and $Y$ be a compact metric space with the topology induced by the metric. Suppose that $f: Y\rightarrow  \iGR{X}$ is continuous. Then the map $f_d\coloneqq \DMod{X}\circ f:Y\rightarrow \iGR{X}$ is a continuous map such that 
	\begin{enumerate}
		\item $Im(f_d)\subset X \subset \iGR{X}$,
		\item If $f(y)\in X \subset \iGR{X}$, then $f_d(y)=f(y)$,
		\item $f_d\simeq f$.
	\end{enumerate}
In addition, if $X$ is locally finite, then $f$ satisfies that
\begin{enumerate}\addtocounter{enumi}{3}
\item $\#(Im(f_d))<\infty$.
\end{enumerate}
\end{lema}

\begin{proof}
First suppose that $X$ is any graph. Let $f:Y\rightarrow \iGR{X}$ be continuous and define $f_d\coloneqq \Theta_{X}\circ f$. Then by definition $Im(f_d)\subset X$ and $f_d(y)=f(y)$ if $f(y)\in X$.

To prove that $f_d\simeq f$, \autoref{prop:SHEBetweenGraphandRGraph} states that $\DMod{X}$ is homotopic to the identity $\Id{\iGR{X}}:\iGR{X}\rightarrow \iGR{X}$. Thus, $f_d = \DMod{X}\circ f \simeq \Id{\iGR{X}}\circ f=f$.

In the final part of the proof, suppose that $X$ is locally finite. Since $f_d$ is continuous, for being the composition of continuous maps, then for every $y\in Y$ there exists $r_y>0$ such that $f_d(\Ball{y}{r_y}{Y})\subset c_{\iGR{X}}(f_d(y))$ because $c_{\iGR{X}}(f_d(y))$ is the minimal neighborhood of $f_d(y)$ in $\iGR{X}$. Since $X$ is locally finite, then $\# f_d(\Ball{y}{r_y}{Y}) \leq \# (c_{\iGR{X}}(f_d(y))\cap X)<\infty$.
Since $Y$ is compact, there exists a finite number of points of $Y$, say $\{y_1,\ldots,y_k\}$, such that $Y=\cup_{i=1}^k \Ball{y_i}{r_{y_i}}{Y}$.
Thus $\# f_d(Y)\leq \sum_{i=1}^k \# f_d(\Ball{y_i}{r_{y_i}}{Y})$, i.e., the set $f_d(Y)$ is finite.
\end{proof}

The next lemma removes the hypothesis of a locally finite graph $X$ from \autoref{lem:DiscretizationOfAMap}, and applies to any graph.
In addition, we consider maps of pairs instead of maps preparing the results for use in homotopy later.
As usual, a map of pairs $f:(Y,A)\to (Z,B)$ is continuous if $f:Y\to Z$ is continuous and $f(A)\subset B$.

\begin{lema}
\label{lem:DiscretizationOfAMapConn}
Let $X$ be a graph and $Y$ be a compact metric space with the topology induced by the metric. Let $A$ be a compact subset of $Y$ and $v\in X$. Suppose that $f:(Y, A)\rightarrow (\iGR{X},v)$ is continuous. Then there exists $f'_d:(Y,A)\rightarrow (\iGR{X},v)$ a continuous map such that
\begin{enumerate}
\item $Im(f'_d)\subset X\subset \iGR{X}$,
\item $\# (Im(f'_d))<\infty$,
\item $f'_d(A)=v$,
\item $f'_d \simeq f$.
\end{enumerate}
\end{lema}

\begin{proof}
Define $f_d\coloneqq \DMod{X}\circ f:Y\rightarrow \iGR{X}$. By \autoref{lem:DiscretizationOfAMap}, we obtain that
\begin{itemize}
\item $Im(f_d)\subset X\subset \iGR{X}$,
\item $f_d(A)=v$, because by hypothesis $f(A)=v$, and
\item $f_d \simeq f$.
\end{itemize}
The following construction is very similar to the subsequent proof of \autoref{lem:FloodPreservingACompact} but with some nuances:
Since for every $x\in X$, $c_{\iGR{X}}(x)$ is the minimal neighborhood of $x$, so for all $y\in Y$ there exists $r_y>0$ such that $f_d(\Ball{y}{r_y}{Y})\subset c_{\iGR{X}}(f_d(y))$.
In addition, since $A$ is a compact subset of a topological Hausdorff space by hypothesis, $A$ is closed, then for every $y\in Y\backslash A$ we can choose $r_y$ sufficiently small so that $\Ball{y}{r_y}{Y}\cap A =\varnothing$.
Then $\{ \Ball{y}{r_y/2}{Y} \mid y\in Y\}$ is a open cover of the compact space $Y$, and there exists a finite number of points $\{y_1,\ldots,y_n\}$ such that
\begin{align*}
Y = \bigcup\limits_{i=1}^n \Ball{y_i}{r_{y_i}/2}{Y}.
\end{align*}
We need two additional claims before we can construct a candidate for $f'_d$.
\begin{enumerate}[wide, labelwidth=!, labelindent=0pt, label=\textbf{Claim \arabic*:}]
	\item If $y_i\notin A$, then $\Ball{y_i}{r_{y_i}/2}{Y}\subset \Ball{y_i}{r_{y_i}}{Y} \subset Y\backslash A$.
	
	\item Consider $y_i$ and $y_j$ such that $\Ball{y_i}{r_{y_i}/2}{Y} \cap \Ball{y_j}{r_{y_j}/2}{Y}\neq \varnothing$, the triangle inequality implies that
	\begin{align*}
		d(y_i,y_j) < (r_{y_i}+r_{y_j})/2 \leq \max(r_{y_i},r_{y_j}),.
	\end{align*}
	and then $y_i \in \Ball{y_j}{r_{y_j}}{Y}$ or $y_j\in \Ball{y_i}{r_{y_i}}{Y}$.
	Both possibilities imply that $f_d(y_i) \in c_{\iGR{X}}(f_d(y_j))$.
\end{enumerate}
The first claim follows directly for the choose of $r_{y_i}$.
The second claim follows because the closure structure on 
$X$ is induced by a graph, and then the closure structure on $X \subset \iGR{X}$ is symmetric, i.e $y_i, y_j \in X \subset \iGR{X}$ and $y_i \in c_{\iGR{X}}(y_j)$ if and only if $y_j \in c_{\iGR{X}}(y_i))$.

For every $i \in \{1,\dots, n\}$, define $f_i:Y \to X \subset \iGR{X}$ by
\begin{align*}
f_i(y)\coloneqq \begin{cases}
f_d(y_i)   & \text{if } y\in \Ball{y_i}{r_{y_i}/2}{Y}, \\
f_{i-1}(y) & \text{otherwise},
\end{cases}
\end{align*}
with $f_0 = f_d$.
Note that, by definition, since $Im(f_d) \subset X \subset \iGR{X}$, it follows that $Im(f_i) \subset X \subset \iGR{X}$ for every $i \in \{1,\dots,n\}$, then $f_n$ satisfies (1).

Since $\{\Ball{y_i}{r_{y_i}/2}{Y}\mid i\in\{1,\ldots,n\}\}$ is a cover of $Y$, then every $y\in Y$ is within a ball of that form.
Thus, $Im(f_n) \subset \{f_d(y_1), \ldots, f_d(y_n) \}$, and $\#Im(f_n) < \infty$.
This paragraph implies that $f_n$ satisfies (2).

We now claim that, for every $i \in \{1,\dots,n\}$, each $f_i$ is continuous, $f_i\simeq f_{i-1}$, and $f_i(A)=v$. 

Claims 1 and 2 above imply that $\Ball{y_i}{r_{y_i}/2}{Y} \subset Y\backslash A$ for any $i \in \{1,\dots,n\}$ with $y_i \in Y\backslash A$.
Therefore for every such $y_i\in Y\backslash A$, $f_i(y) = f_{i-1}(y)$ for any point $y \in A$.
On the other hand, $f_i(y_i)=v$ for any $y_i\in A$.
Taking into consideration that $A$ is covered by the $\Ball{y_i}{r_{y_i}}{Y}$ with $y_i \in A$ by construction, we obtain that $f_i(A) = v$.
These two facts imply that $f_i(A) = v$ for every $i\in \{1,\ldots,n\}$. In particular, $f_n$ satisfies (3).

Now, by induction on $i$, we prove that $f_i$ is continuous. First note that $f_0$ is continuous by \autoref{lem:DiscretizationOfAMap}, thus the base case holds. For the inductive step, suppose that $f_{i-1}$ is continuous.

Recalling that in a metric space, a filter base for the neighborhood filter $\NFil{y}$ is the set of all the $\Ball{y}{r}{Y}$ with $r>0$, we will prove the continuity of $f_i$ considering the following three cases over $y\in Y$:
\begin{enumerate}[wide, labelwidth=!, labelindent=0pt, label=\textbf{Case \arabic*:}]

\item Let $y\in \Ball{y_i}{r_{y_i}/2}{Y}$. For $Y$ being a metric space, there exists $r>0$ such that $\Ball{y}{r}{Y}\subset \Ball{y_i}{r_{y_i}/2}{Y}$, and for the definition of $f_i$, we have that $f_i(z)=f_i(y_i)=f_d(y_i)$ for every $z\in \Ball{y}{r}{Y}$. Thus, \[f_i(\Ball{y}{r}{Y}) = f_i(\Ball{y_i}{r_{y_i}/2}{Y}) = \{f_i(y_i)\}  \subset c_{\iGR{X}}(f_i(y_i))=c_{\iGR{X}}(f_i(y)).\]
Concluding that $f_i(\NFil{y}) = \GFil{f_i(y)}{\uGRCC{X}} \supset \GFil{c_{\iGR{X}}(f_i(y))}{\uGRCC{X}}$ and thus $f_i$ is continuous at $y$.

\item Let $y\in \partial \Ball{y_i}{r_{y_i}/2}{Y}$. By the definition of the boundary of a set, there exists $r>0$ such that $\Ball{y}{r}{Y} \subset \Ball{y_i}{r_{y_i}}{Y}$. Since $f_i$ is continuous by the inductive hypothesis, we can also ask for $r$ satisfying
$f_{i-1}(\Ball{y}{r}{Y}) \subset c_{\iGR{X}}(f_{i-1}(y))$. 

Now define the set $\Lambda\coloneqq \{k\mid y\in \Ball{y_k}{r_{y_k}}{Y}\} \subset \{1,\ldots,n\}$ and note that $f_{i-1}(y)\in \{f_d(y)\}\cup \{ f_d(y_k)\mid k\in \Lambda\}$ and $y \in \Ball{y_k}{r_{y_k}}{Y} \cap \Ball{y_i}{r_{y_i}}{Y} \neq \varnothing$ for every $k\in \Lambda$. 

In addition, $y\in \Ball{y_i}{r_{y_i}}{Y}$, and therefore $f_d(y)\in c_{\iGR{X}}(f_d(y_i))$ by the choice of $r_{y_i}$. Thus, given that the closure structure on $X\subset \iGR{X}$ is symmetric, it follows that $f_d(y_i)\in c_{\iGR{X}}(f_{i-1}(y))$. Putting this together, and noting that $f_{i-1}(y) = f_i(y)$ because $y\notin \Ball{y_i}{r_{y_i}/2}{Y}$, we have
\begin{align*}
f_i(\Ball{y}{r}{Y}) \subset & f_{i-1}(\Ball{y}{r}{Y}\cup \{f_d(y_i)\}) \subset c_{\iGR{X}}(f_{i-1}(y))\cup \{f_d(y_i)\} \subset \\ & c_{\iGR{X}}(f_{i-1}(y)) = c_{\iGR{X}}(f_i(y)).
\end{align*}
Thus $f_i$ is continuous ay $y$ for the same argument in the end of \textbf{Case 1}.

\item Let $y \notin \overline{\Ball{y_i}{r_{y_i}/2}{Y}}$. Then $\Ball{y}{r}{Y}\subset Y\backslash\overline{\Ball{y_i}{r_{y_i}/2}{Y}}$ for some $r>0$. Thus \[f_i(\Ball{y}{r}{Y})=f_{i-1}(\Ball{y}{r}{Y}) \subset c_{\iGR{X}}(f_{i-1}(y)) = c_{\iGR{X}}(f_i(y)).\] Also concluding that $f_i$ is continuous at $y$ for the argument in the previous two cases.
\end{enumerate}
The three possible cases implies that $f_{i}$ is continuous if $f_{i-1}$. In particular, the map $f_n:Y\to \iGR{X}$ is continuous.

Finally, we prove there is a homotopy between $f_{i-1}$ and $f_{i}$. To achieve this, let's define a map $H_i:Y\times I\rightarrow X\subset \iGR{X}$ such that
\begin{align*}
H_i(y,t) \coloneqq \left\lbrace\begin{array}{ll}
f_i(y)     & t=1, \\
f_{i-1}(y) & t\in [0,1),
\end{array}\right.
\end{align*}
for $i\in\{1,\ldots,n\}$. Note that for $t\in [0,1)$ we can find $r>0$ small enough such $H_i(\Ball{(y,t)}{r}{Y\times I})=f_{i-1}(\Ball{(y,t)}{r}{Y\times I})$, and therefore
$H$ is continuous at any point $(y,t)$ in $Y \times [0,1)$.

We now show that $H$ is continuous at any point of the form $(y,1) \in Y \times [0,1]$, considering the same three cases as above:
\begin{enumerate}[wide, labelwidth=!, labelindent=0pt, label=\textbf{Case \arabic*:}]
\item Let $y \in\Ball{y_i}{r_{y_i}/2}{Y}$. Define $\Lambda\coloneqq \{k\mid y\in \Ball{y_k}{r_{y_k}/2}{Y}\}$, then 
\begin{equation*}
f_d(y) \in f(\Ball{y_k}{r_{y_k}/2}{Y}) \subset c_{\iGR{X}}(f_d(y_k))
\end{equation*}
 for every $k\in\Lambda$. Since $f_d$ is continuous, there is a $r>0$ such that 
 $f_d(\Ball{y}{r}{Y})\subset c_{\iGR{X}}(f_d(y))$. In addition, we can choose $r$ such that $\Ball{y}{r}{Y} \subset \Ball{y_k}{r_{y_k}/2}{Y}$ for every $k\in\Lambda$ given that $Y$ is a metric space.  Recall that, $\forall y \in Y$, $r_y>0$ was chosen so that $f_d(\Ball{y}{r_y}{Y}) \subset c_{\iGR{X}}(f_d(y))$. As a consequence, $\Ball{y}{r}{Y} \subset \Ball{y_i}{r_{y_i}/2}{Y}$ implies that $f_d(\Ball{y}{r}{Y}) \subset c_{\iGR{X}}(f_d(y_i))$. Also, by construction, $\Ball{y_k}{r_{y_k}}{Y} \cap \Ball{y_i}{r_{y_i}}{Y} \neq \varnothing$
 for all $k \in \Lambda$, and from the remarks above, it follows that $f_d(y_k) \in c_{\iGR{X}}(f_d(y_i))$ for all $k\in \Lambda$. This implies that 
 \begin{align*}
 f_{i-1}(\Ball{y}{r}{Y})\subset &f_d(\Ball{y}{r}{Y})\cup \{f_d(y_j)\mid
  j\in\Lambda\}\subset\\ & c_{\iGR{X}}(f_d(y_i)) \cup \{f_d(y_j)\mid j\in\Lambda\}\subset c_{\iGR{X}}(f_d(y_i)).
  \end{align*} 
Furthermore, $f_i(y) = f_d(y_i)$ and $f_i(\Ball{y}{r}{Y}) = \{f_d(y_i)\}$ by the definition of $f_i$ and the choice of the radius $r$. We conclude that 
\begin{align*}
H_i(\Ball{(y,1)}{r}{Y\times I})\subset & f_i(\Ball{y}{r}{Y})\cup f_{i-1}(\Ball{y}{r}{Y})\subset c_{\iGR{X}}(f_d(y_i)) = \\
& c_{\iGR{X}}(f_i(y)) = c_{\iGR{X}}(H_i(y,1)).
\end{align*}
Hence $H$ is continuous at $(y,1)$.
\item Let $y \in \partial \Ball{y_i}{r_{y_i}/2}{Y}$.
As a consequence of $f_i$ and $f_{i-1}$ being continuous,  there is a $r>0$ such that $f_i(\Ball{y}{r}{Y})\subset c_{\iGR{X}}(f_i(y))$ and $f_{i-1}(\Ball{y}{r}{Y})\subset c_{\iGR{X}}(f_{i-1}(y)) = c_{\iGR{X}}(f_{i}(y))$.
Concluding \[H_i(\Ball{(y,1)}{r}{Y\times I})\subset f_i(\Ball{y}{r}{Y})\cup f_{i-1}(\Ball{y}{r}{Y})\subset c_{\iGR{X}}(f_i(y)),\]
and hence $H$ is continuous at $(y,1)$.

\item Let $y \notin \overline{\Ball{y_i}{r_{y_i}/2}{Y}}$. Since the closure of the ball is a closed set, then there is $r>0$ such that $y$ does not meet $\Ball{y_i}{r_{y_i}/2}{Y}$, and thus
\begin{align*}
H_i(\Ball{(y,1)}{r}{Y\times I}) =& f_{i}(\Ball{y}{r}{Y})=f_{i-1}(\Ball{y}{r}{Y}) \subset c_{\iGR{X}}(f_{i-1}(\Ball{y}{r}{Y}) =\\ & c_{\iGR{X}}(f_i(\Ball{y}{r}{Y}) = c_{\iGR{X}}(H_i(\Ball{(y,1)}{r}{Y\times I})).
\end{align*}
Also concluding that $H$ is continuous at $(y,1)$.
\end{enumerate}

The three possible cases imply that $H_i$ is continuous for every $i\in \{1,\ldots,n\}$. This means that $f_i\simeq f_{i-1}$. Therefore, for the transitivity of homotopy relation, $f_n\simeq f_d \simeq f$. Consequently, $f_n$ satisfies (4).

Defining $f'_d \coloneqq f_n$, we find the desired map.
\end{proof}

\begin{definicion}\label{def:DiscModMap}
	Let $X$ be a graph and $Y$ be a compact metric space with the topology induced by the metric. Let $A$ be compact subset of $Y$ and $v\in X$. Suppose that $f:(Y, A)\rightarrow (\iGR{X},v)$ is continuous. The map $f'_d$ from \autoref{lem:DiscretizationOfAMapConn} is said to be \emph{the finite discrete modification of $f$}, denoted by $\DMod{X}[f]$.\nomenclature{$\DMod{X}[f]$}{The discrete modification of a map $f$}
\end{definicion}

	Note that the finite discrete modification converts an infinite graph problem into a finite graph problem, which is a subspace of the original graph.
	Additionally, the discrete modification is uniquely defined consider our graphs with well-ordered vertices.
	
	Even though the graphs used in the rest of the document are not necessarily locally finite, and even though \autoref{lem:DiscretizationOfAMapConn} generalizes that case, we decided to leave that part of \autoref{lem:DiscretizationOfAMap} as a vestige of our initial approach to the main theorem of this article.

	The maps $f_i$ in the proof of \autoref{lem:DiscretizationOfAMapConn} are examples of the following definition:
	The $v$-flood.
	This map is an adaptation from the $v$-flood in digital topology \cite{Lupton2023DigitalHom}.

\begin{definicion}[The $v$-flood]
\label{def:Flood}
	Given a continuous map $f:Y\rightarrow \iGR{X}$ and $v\in X \subset \iGR{X}$, and let $B\subset Y$.
	We define the not-necessarily continuous function $\Fld{f}{B}{v}:Y\rightarrow \uGRCC{X}$
\begin{align*}
\Fld{f}{B}{v}(y)\coloneqq \left\lbrace \begin{array}{ll}
v & y\in B,\\
f(y) & \text{otherwise},
\end{array}\right.
\end{align*}
which we call \emph{the flood of $f$ over $B$ with $v$} or \emph{the $v$-flood of $f$ over $B$}.\nomenclature{$\Fld{f}{B}{v}$}{The $v$-flood of $f$ over $B$}
\end{definicion}

\begin{figure}[h]
\begin{tikzpicture}
\fill[NavyBlue] (-1,1) -- (1,1) -- (1,0.3) -- (-1, -0.3);
\fill[Green] (-1,-1) -- (1,-1) -- (1,0.3) -- (-1, -0.3);
\draw[Red, very thick] (1,0.3) -- (-1,-0.3);

\draw[->] (1.3,0)--(1.7,0);

\fill[NavyBlue] (2,1) -- (4,1) -- (4,0.3) -- (2, -0.3);
\fill[Green] (2,-1) -- (4,-1) -- (4,0.3) -- (2, -0.3);
\draw[Red, very thick] (4,0.3) -- (2,-0.3);
\fill[Black, semitransparent] (2,-0.6) -- (2,0) -- (4,0.6) -- (4,0);

\draw[->] (4.3,0)--(4.7,0);

\fill[NavyBlue] (5,1) -- (7,1) -- (7,0.3) -- (5, -0.3);
\fill[Green] (5,-1) -- (7,-1) -- (7,0.3) -- (5, -0.3);
\fill[Red] (5,-0.6) -- (5,0) -- (7,0.6) -- (7,0);

\draw[black] (-1,-1) rectangle (1,1);
\draw[black] (2,-1) rectangle (4,1);
\draw[black] (5,-1) rectangle (7,1);
\end{tikzpicture}
\caption{Illustration of a flood.}
\label{fig:Flood}
\end{figure}
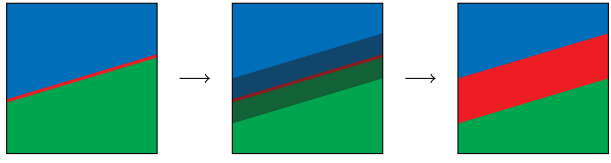

\autoref{fig:Flood} shows subset of $Y$. The colors blue, green and red representing points in $X$.
In particular, $v$ is red.
The first square represents the image of $f$, and the gray transparent zone in the second square is the set $B$ in that subset of $Y$.
The last picture represents the image of $\Fld{f}{B}{v}$.

\begin{lema}
\label{lem:Flood}
Let $Y$ be a metric space with the topology induced by the metric. Let $f: Y\rightarrow \iGR{X}$ be a continuous map such that $Im(f)\subset X \subset \iGR{X}$ and let $v\in X \subset \iGR{X}$. If for every $y\in Y$ such that $f(y)=v$, there exists a $r_y>0$ such that $f(\Ball{y}{r_y}{Y})\subset c_{\iGR{X}}(v)$, then $\Fld{f}{B}{v}$ is continuous in $X\subset \iGR{X}$ and $f\simeq \Fld{f}{B}{v}$ for
\begin{align*}
B \coloneqq \bigcup_{y\in f^{-1}(v)} \Ball{y}{r_y/2}{Y}.
\end{align*}
And thus it follows the same result replacing in $B$ the elements $r_y$ for a fixed $r'_y\leq r_y$.
\end{lema}

\begin{proof}
Suppose that for every $y\in Y$ such that $f(y)=v$ there exists a $r_y>0$ such that $f(\Ball{y}{r_y}{Y})\subset c_{\iGR{X}}(v)$.
First we prove that $\Fld{f}{B}{v}$ is continuous. We consider two cases over $y\in Y$:

\begin{enumerate}[wide, labelwidth=!, labelindent=0pt, label=\textbf{Case \arabic*:}]
\item Let $y\in B$. Since $B$ is the union of open sets, then $B$ is open in $Y$, and there exists $r>0$, in particular we can choose $r<r_y$, such that $\Ball{y}{r}{Y} \subset B$ and therefore $\Fld{f}{B}{v}(\Ball{y}{r}{Y})=v \in c_{\iGR{X}}(v)$. Hence  $\Fld{f}{B}{v}(\Ball{y}{r}{Y}) \subset c_{\iGR{X}}(v)$, and $\Fld{f}{B}{v}$ is continuous at $y$.

\item Let $y\in Y\backslash B$. By the continuity of $f$, and being $c_{\iGR{X}}(f(y))$ the minimal neighborhood of $f(y)$ in $\iGR{X}$, $f(\Ball{y}{r}{Y})\subset c_{\iGR{X}}(f(y))$ for some $r>0$. Let's consider the two possible subcases for $v$:
\begin{enumerate}
\item Let $v\in c_{\iGR{X}}(f(y))$. Then \[\Fld{f}{B}{v}(\Ball{y}{r}{Y})\subset f(\Ball{y}{r}{Y})\cup\{v\}\subset c_{\iGR{X}}(f(y)).\]

\item Let $v\notin c_{\iGR{X}}(f(y))$. Since $c_{\iGR{X}}(f(y))$ is the minimal neighborhood, there exists $r'>0$ such that $f(\Ball{y}{r'}{Y})$ does not contain the elements which are not neighbors of $f(y)$; in particular, we can choose $r'<r$. Then $v\notin f(\Ball{y}{r'}{Y})$. That $r'$ has to be less or equal to the distance from $y$ to any point with image $v$ under $f$; thus $\Ball{y}{r'}{Y}\subset Y\backslash B$ and $\Fld{f}{B}{v}(z) = f(z)$ for any $z\in \Ball{y}{r'}{Y}$.
Therefore
\[\Fld{f}{B}{v}(\Ball{y}{r'}{Y}) =  f(\Ball{y}{r'}{Y})\subset  c_{\iGR{X}}(f(y))=c_{\iGR{X}}(\Fld{f}{B}{v}(y)).\]

\end{enumerate}
In subcases (a) and (b), we conclude that $\Fld{f}{B}{v}$ is continuous at $y$.
\end{enumerate}
Combining the two cases over $y\in Y$, we conclude that $\Fld{f}{B}{v}$ is continuous.

The next step is to define the homotopy. We define $H:Y\times I\rightarrow X\subset \iGR{X}$ by
\begin{align*}
H(y,t) \coloneqq \left\lbrace\begin{array}{ll}
f(y)   & \text{ if }t\in [0,1),\\
\Fld{f}{B}{v}(y) & \text{ if }t=1.
\end{array}\right.
\end{align*}
By similar arguments to the ones in \autoref{prop:SHEBetweenGraphandRGraph}, note that $H$ is continuous at any point in $Y\times [0,1)$. Now, for $t=1$, we consider two cases over $f(y)$:

\begin{enumerate}[wide, labelwidth=!, labelindent=0pt, label=\textbf{Case \arabic*:}]
\item Let $f(y)=\Fld{f}{B}{v}(y)$. By the continuity of $f$, there is a $r>0$ such that $f(\Ball{y}{r}{Y})\subset c_{\iGR{X}}(f(y))=c_{\iGR{X}}(\Fld{f}{B}{v}(y))$. Further, for the continuity of $\Fld{f}{B}{v}$, we can choose $r$ small enough such that $\Fld{f}{B}{v}(\Ball{y}{r}{Y})\subset c_{\iGR{X}}(\Fld{f}{B}{v}(y))$.
Thus
\begin{align*}
	H(\Ball{(y,1)}{r}{Y\times I})\subset & f(\Ball{y}{r}{Y})\cup \Fld{f}{B}{v}(\Ball{y}{r}{Y})\subset \\ & c_{\iGR{X}}(\Fld{f}{B}{v}(y))=c_{\iGR{X}}(H(y,1)).
\end{align*}
Hence $H$ is continuous at $(y,1)$.

\item Let $f(y)\neq \Fld{f}{B}{v}(y)$. Given the definition of $\Fld{f}{B}{v}$, $y$ must be an element of $B$. This implies that there is a $y_0\in Y$ such that $f(y_0) = v$ and $y\in \Ball{y_0}{r_{y_0}/2}{Y}$.
Since $\Ball{y_0}{r_{y_0}/2}{Y}$ is open, there is $r$ such that $\Ball{y}{r}{Y} \subset \Ball{y_0}{r_{y_0}/2}{Y}$.

This implies that $\Fld{f}{B}{v}(\Ball{y}{r}{Y})=\{v\}$, for the definition of the $v$-flood over $B$, and $f(\Ball{y}{r}{Y}) \subset c_{\iGR{X}}(v)$ for the choose of $r_{y_0}$. Then
\begin{align*}
	H(\Ball{(y,1)}{r}{Y\times I})\subset & f(\Ball{y}{r}{Y})\cup \{ v\} = c_{\iGR{X}}(v) = \\ & c_{\iGR{X}}(\Fld{f}{B}{v}(y))=c_{\iGR{X}}(H(y,1)).
\end{align*}
Hence $H$ is continuous at $(y,1)$.
\end{enumerate} 
Thus, combining the two cases, $H$ is continuous. Concluding that $f \simeq \Fld{f}{B}{v}$.
\end{proof}

Now, we will prove an extension of \autoref{lem:Flood} for maps of pairs. For a domain $(Y,A)$, this version guarantees that the image of $(y,t)$ remains constant for all $t$.

\begin{lema}
\label{lem:FloodPreservingACompact}
Let $Y$ be a metric space with topology induced by the metric. Let $A$ be a compact subset of $Y$ and $v,w\in X\subset \iGR{X}$. Suppose $f:(Y,A)\rightarrow (\iGR{X},w)$ a continuous function such that $Im(f)\subset X \subset \iGR{X}$. If for every $y\in Y$ such that $f(y)=v$, there exists $r_y>0$ such that
\begin{itemize}
\item $f(\Ball{y}{r_y}{Y})\subset c_{\iGR{X}}(v)$ and
\item if $v\neq w$, $A\cap \Ball{y}{r_y}{Y}=\varnothing$,
\end{itemize}
then $\Fld{f}{B}{v}(A)=w$, $\Fld{f}{B}{v}$ is continuous in $X\subset \iGR{X}$ and $f\simeq \Fld{f}{B}{v} \operatorname{rel}(A)$ in the closure space $\iGR{X}$ for
\begin{align*}
B \coloneqq \bigcup_{y\in f^{-1}(v)} \Ball{y}{r'_y/2}{Y}
\end{align*}
with $r'_y\leq r_y$.
\end{lema}

\begin{proof}
First we prove that $\Fld{f}{B}{v}(A)=w$.
If $v=w$ there is nothing to do, by hypothesis $\Fld{f}{B}{v}(f^{-1}(v))=\{v\}$ and $A\subset f^{-1}(v)$.
On the other hand, if $v\neq w$, then $B\subset Y\backslash A$ by the choose of the $r_y$.
Thus $f^v_B(A)=f(A)=w$ by definition.

For \autoref{lem:Flood}, $\Fld{f}{B}{v}$ is already continuous and $f\simeq \Fld{f}{B}{v}$. Considering the map $H$ defined in \autoref{lem:Flood}, $H(y,1)=\Fld{f}{B}{v}(y)$ and $H(y,t)=f(y)$ for every $t\in [0,1)$. Thus $H(A,t)=w$, and then the homotopy is relative to $A$.
\end{proof}

\begin{observacion}\label{obs:BuildFldDisMod}
	Given a compact metric space $Y$ and $A\subset Y$ compact. For any pair map $f: (Y,A)\rightarrow (\iGR{X},w)$ with $w\in X$, by \autoref{lem:DiscretizationOfAMapConn} we have that  $Im(\DMod{X}[f]) = \{v_0,\ldots,v_n\} \subset X$ is finite and $\DMod{X}[f](A)=w$. In addition, since $f$ is continuous and $A$ is compact, for every $y\in Y$ such that $f(y)= v_0$, we can find $r_y>0$ such that $f(\Ball{y}{r_y}{Y})\subset c_{\iGR{X}}(v_0)$ and such that $A\cap \Ball{y}{r_y}{Y}=\varnothing$ if $v_0\neq w$.
	Thus, by \autoref{lem:FloodPreservingACompact}, we conclude that $\Fld{\DMod{X}[f]}{B_0}{v_0} \simeq \DMod{X}[f]\simeq f$ and $\Fld{\DMod{X}[f]}{B_0}{v_0}(A)=w$ with
\begin{align*}
B_0\coloneqq\bigcup_{y\in (\DMod{X}[f])^{-1}(v_0)} \Ball{y}{r_y/2}{Y}.
\end{align*}
Note we can repeat this process considering the map $\Fld{\DMod{X}[f]}{B_0}{v_0}$, the point $v_1$ and building $B_1\subset Y$ in a similar fashion as $B_0$; obtaining $\Fld{\Fld{\DMod{X}[f]}{B_0}{v_0}}{B_1}{v_1}$ which is homotopic to $\Fld{\DMod{X}[f]}{B_0}{v_0}$ and preserves the image of $A$.

Then we inductively repeat this process a finite number of times for all of the vertices in $Im(\DMod{X}[f])$, obtaining
\begin{align*}
f \simeq & \Fld{\DMod{X}[f]}{B_0}{v_0} \simeq \Fld{\Fld{\DMod{X}[f]}{B_0}{v_0}}{B_1}{v_1} \simeq \ldots \simeq \Fld{\cdots \Fld{\Fld{\DMod{X}[f]}{B_0}{v_0}}{B_1}{v_1} \cdots}{B_n}{v_n},\\
\{w\} = & \Fld{\cdots \Fld{\Fld{\DMod{X}[f]}{B_0}{v_0}}{B_1}{v_1} \cdots}{B_n}{v_n}(A).
\end{align*}
\end{observacion}

\begin{definicion}\label{def:FldDiscMod}
	Let $Y$, $A$, $X$ and $f$ as in \autoref{obs:BuildFldDisMod} such that the image of $\DMod{X}[f]$ is the finite set $\{v_0,\ldots,v_n\}$. For $m\in\{0,\ldots,n\}$, the map \[\FldF{f}{v_m} = \Fld{\cdots \Fld{\Fld{\DMod{X}[f]}{B_0}{v_0}}{B_1}{v_1} \cdots}{B_m}{v_m} \] is said to be the \emph{flooding discrete modification of $f$ until $v_m$}.\nomenclature{$\FldF{f}{v_m}$}{The flooding discrete modification of $f$ until $v_m$}
\end{definicion}

In order to have the document auto-contained or referenced, we add the following small result in the context of metric spaces.

\begin{proposicion}
\label{prop:EventuallyConstant}
Suppose that we have the sequences $\{y_n\}\rightarrow x$ and $\{z_n\}\rightarrow x$ in a metric space $Y$ and we have a function $f:Y\rightarrow X$ (not necessarily continuous) such that $\#X<\infty$. Then we can take a subsequence of $\{(y_n,z_n)\}$ such that $f(y_n)=f(y_k)$ and $f(z_n)=f(z_k)$ for every $k,n$.
\end{proposicion}

\begin{proof}
Note that the image of the sequence contains at most $(\#X)^2$ points. If we suppose that we cannot make that construction, then for every $(v,w)\in X\times X$ there exists a $r_{v,w}$ such that through that point the image cannot be $(v,w)$. If we take $r\coloneqq \min_{v,w}r_{v,w}$ $\Absurd$. From that $r$, we cannot take more points and then their images are out of $X$.
\end{proof}

We use the previous result just in the following lemma.

\begin{lema}
\label{lem:SimplexIsInsideANeighborhood}
Let $Y$ be a compact metric space and $X$ a graph. Given a continuous map $f:Y\rightarrow \iGR{X}$ with image $\{v_0,\ldots,v_m\}$, for every $x\in Y$ there exists $r_x>0$ such that if $v,w\in \FldF{f}{v_m}(\Ball{x}{r_x}{Y})$, then $v\in c_{\iGR{X}}(w)$ (equivalently $w\in c_{\iGR{X}}(v)$).
\end{lema}

\begin{proof}
We proceed by contradiction.
Suppose that there exists $x\in Y$ such that for every $r>0$ there exists a pair $w_r,w'_r\in \FldF{f}{v_m}(\Ball{x}{r}{Y})$ we have that $w_r\notin c_{\iGR{X}}(w'_r)$.
We can restrict this sentence for the numbers $1/k$ for every $k\in\mathbb{N}$.
For this reason, we have elements $y_k,z_k\in X$ such that $y_k,z_k\in \Ball{x}{\tfrac{1}{k}}{Y}$, $\FldF{f}{v_m}(y_k)=w_{1/k}$ and $\FldF{f}{v_m}(z_k)=w'_{1/k}$.
Thus we have two sequences $\{y_k\}$ and $\{z_k\}$ that converge to $x$.
By \autoref{prop:EventuallyConstant}, we can take a subsequence of both such that the image of all of the elements of the subsequence $\{y_k\}$ (resp. $\{z_k\}$) is a fixed element $w$ (resp. $w'$).

We name $v\coloneqq \FldF{f}{v_m}(x)$, and divide this proof in two cases:
\begin{enumerate}[wide, labelwidth=!, labelindent=0pt, label=\textbf{Case \arabic*:}]
\item Let $w<v$ or $w'<v$ $\Absurd$.  We obtain a direct contradiction because the elements only change to $w$ (resp. $w'$) in the flooding $\FldF{f}{w}$ (resp. $\FldF{f}{w'}$). Since $w<v$ (or $w'<v$), $\FldF{f}{v}(\Ball{x}{r}{Y})=\{v\}$ and they will not change to $w$ (or $w'$) after this point.

\autoref{fig:vGreaterThanw} illustrates this first case.
Each rectangle represents the image of a subset of $Y$ with the colors blue, green and red representing $w$, $w'$ and $v$ respectively.
The white point with black contour is $x$. From left to right: The first rectangle is the image before $\FldF{f}{v}$. The second rectangle adds a point whose image is $v$ and its neighborhood to the first rectangle.
The last rectangle shows $\FldF{f}{v}$ and illustrates that there is a neighborhood around $x$ without blue or green points.
For illustrative purposes, all the other points in the image were ignored.

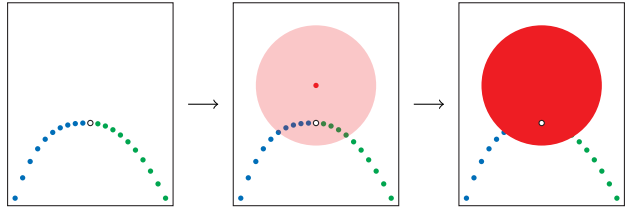
\begin{figure}[h]
\begin{tikzpicture}
\foreach \t in {1,...,10}
	\fill[Green] (\t/10,{-(\t/10)^2}) circle (1pt);
\foreach \t in {-10,...,-1}
	\fill[NavyBlue] (\t/10,{(\t/10)^3}) circle (1pt);
\fill[white] (0,0) circle (1pt);
\draw[Black] (0,0) circle (1pt);

\draw[->] (1.3,.25)--(1.7,.25);

\foreach \t in {1,...,10}
	\fill[Green] (\t/10+3,{-(\t/10)^2}) circle (1pt);
\foreach \t in {-10,...,-1}
	\fill[NavyBlue] (\t/10+3,{(\t/10)^3}) circle (1pt);
\fill[Red] (3,1/2) circle (1pt);
\fill[Red, nearly transparent] (3,1/2) circle (.8cm);
\fill[white] (3,0) circle (1pt);
\draw[Black] (3,0) circle (1pt);

\draw[->] (4.3,.25)--(4.7,.25);

\foreach \t in {1,...,10}
	\fill[Green] (\t/10+6,{-(\t/10)^2}) circle (1pt);
\foreach \t in {-10,...,-1}
	\fill[NavyBlue] (\t/10+6,{(\t/10)^3}) circle (1pt);
\fill[Red] (6,1/2) circle (.8cm);
\fill[white] (6,0) circle (1pt);
\draw[Black] (6,0) circle (1pt);

\draw[Black] (-1.1,-1.1) rectangle (1.1,1.6);
\draw[Black] (1.9,-1.1) rectangle (4.1,1.6);
\draw[Black] (4.9,-1.1) rectangle (7.1,1.6);
\end{tikzpicture}
\caption{Illustration of $w<v$ or $w'<v$.}
\label{fig:vGreaterThanw}
\end{figure}

\item Suppose that $w>v$ and $w'>v$, without loss of generality, assume that $w<w'$.
Note that in the map $\FldF{f}{v}$, there is an $r>0$ such that the image of $\Ball{x}{r}{Y}$ is $v$; therefore, in $\FldF{f}{w}$, the image of $\{y_k\}$ must become $w$.
Take any element $z\in Y$ such that $\FldF{f}{w}(z)=w'$, then
\[\lim_{n}d(y_k,z) \rightarrow d(x,z),\]
by the continuity of $d$ in the metric space.
For the construction of $\FldF{f}{w'}$, we see that the image of the points outside of $\Ball{z}{d(x,z)/2}{Y}$ does not change to $w'$ in $\FldF{f}{w'}$.
This implies that there are no points inside $\Ball{x}{r/2}{Y}$ with image $w'$ in $\FldF{f}{w'}$, and therefore bot in $\FldF{f}{v_m}$ neither. $\Absurd$

To illustrate this second case, see \autoref{fig:wGreaterThanv}.
Each rectangle represents the image of a subset of $Y$ with the colors blue, green and red representing $w$, $w'$ and $v$ respectively. 
The white point with a black contour is $x$.
From left to right: The first rectangle represents $\FldF{f}{w}$, for which we require a sequence with image $w$.
We assume that there is a green point in the contour of the red circle.
The second rectangle has two circles that illustrate the distance and the half of the distance between the green point in the contour and the white point $x$.
The last rectangle represents $\FldF{f}{w'}$ and shows the space without green points and with blue points.
For illustrative purposes, all the other points in the image were ignored.

\begin{figure}[h]
\begin{tikzpicture}
\fill[Red] (0,0.3) circle (.8cm);
\foreach \t in {-10,...,-1}
	\fill[NavyBlue] (\t/10,{-(\t/10)^3}) circle (1pt);
\fill[white] (0,0) circle (1pt);
\draw[Black] (0,0) circle (1pt);
\fill[Green] (0,-0.5) circle (1pt);

\draw[->] (1.3,.25)--(1.7,.25);

\fill[Red] (3,0.3) circle (.8cm);
\foreach \t in {-10,...,-1}
	\fill[NavyBlue] (\t/10+3,{-(\t/10)^3}) circle (1pt);
\draw[Green,dashed,thick] (3,-0.5) circle (0.5cm);
\draw[Green,dashed,thick] (3,-0.5) circle (0.25cm);
\fill[white] (3,0) circle (1pt);
\draw[Black] (3,0) circle (1pt);
\fill[Green] (3,-0.5) circle (1pt);

\draw[->] (4.3,.25)--(4.7,.25);

\fill[Red] (6,0.3) circle (.8cm);
\foreach \t in {-10,...,-1}
	\fill[NavyBlue] (\t/10+6,{-(\t/10)^3}) circle (1pt);
\fill[Green] (6,-0.5) circle (0.25cm);
\fill[white] (6,0) circle (1pt);
\draw[Black] (6,0) circle (1pt);
\draw[white,dashed,thick] (6,0) circle (0.25cm);

\draw[Black] (-1.1,-1.1) rectangle (1.1,1.6);
\draw[Black] (1.9,-1.1) rectangle (4.1,1.6);
\draw[Black] (4.9,-1.1) rectangle (7.1,1.6);
\end{tikzpicture}
\caption{Illustration of $v<w<w'$.}
\label{fig:wGreaterThanv}
\end{figure}
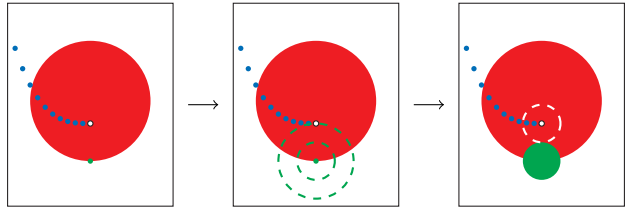
\end{enumerate}

Note that $v=w$ or $v=w'$ is an automatic contradiction.
We conclude that there exist the desired $r_x$.
\end{proof}

We build the discrete modification (\autoref{def:DiscModification}) to convert all of the points in $\uGRCC{X}$ into vertices in $X$.
The next definition describes an almost inverse process, moving from a function whose image is only in $X$ to a map whose image is at $\uGRCC{X}$.

\begin{definicion}[The Convex Transformation]
\label{def:ConvexTransformation}
Let $X$ be a graph, be $Y$ a triangularizable compact space with a finite triangulation $K$, let $\mathcal{U}$ be an open cover of $Y$ and let $f: Y\rightarrow X\subset \iGR{X}$ such that:
\begin{enumerate}
\item For every $U\in\mathcal{U}$, all $y,y'\in U$ satisfies that $f(y)\in c_X(f(y'))$.
\item Every simplex in $Y$ is totally contained in some $U\in\mathcal{U}$.
\end{enumerate}
We define the \emph{convex transformation of $f$ over $K$}, denoted by $\ConM{f}{K}:Y\rightarrow \GRCC{X}$, such that for every simplex $\sigma$ in the triangulation of $Y$, say $\{y_0,\ldots,y_n\}$,
\nomenclature{$\ConM{f}{K}$}{The convex transformation of the map $f$ over the triangulation $K$}
\begin{align*}
\ResF{\ConM{f}{K}}{|\sigma|}\left(\sum_{i=0}^n t_iy_i\right) = & \sum_{i=0}^n t_if(y_i) \\
\text{such that } & t_i\geq 0 \text{ and } \sum_{i=0}^n t_i=1
\end{align*}
the convex sums $\sum t_iy_i$ and $\sum t_if(y_i)$ are an abuse of notation for convenience; $\sum t_iy_i$ means the image of the point in the triangulation in $Y$, and $\sum t_if(y_i)$ means the convex sum of the vertices $f(y_i)$ seen in $\GRCC{X}$, which has the same underlying set that $\iGR{X}$.

When the context is clear, we omit the triangulation and write $\ConM{f}{}$.
\end{definicion}

\begin{lema}
\label{lem:WellDefinedConvTrans}
Consider $X$ a graph, $Y$ a triangularizable compact space with a finite triangulation $K$, $\mathcal{U}$ and $f$ as in \autoref{def:ConvexTransformation}. The map $\ConM{f}{K}:Y\rightarrow \GRCC{X}$ is well-defined.
\end{lema}

\begin{proof}
To prove that $\ConM{f}{K}$ is well-defined, it is enough to show that if two simplexes, say $\Delta^n=\{v_0,\ldots,v_n\}$ and $\Delta^m=\{w_0,\ldots,w_m\}$, share a face, for example $\Delta^\ell=\{x_0,\ldots,x_\ell\}$ such that $x_i=v_i=w_i$ for $i\in [\ell]$ (after a possible relabeling), then the value $\ResF{\ConM{f}{K}}{|\Delta^n|}(y) = \ResF{\ConM{f}{K}}{|\Delta^m|}(y)$ for every $y\in |\Delta^\ell|$.
Let $t_i,s_j,r_k\in [0,1]$ for every $i\in [n],j\in [m],k\in [\ell]$ such that $\sum_{i=0}^n t_i = \sum_{j=0}^m s_j = \sum_{k=0}^\ell r_k=1$ and $t_i=s_j=0$ for $i\in [n]\backslash[\ell],j\in [n]\backslash[\ell]$ we that
\begin{align*}
\ConM{f}{K}\mid_{\Delta^\ell\subset \Delta^n}\left(\sum_{i=0}^n t_ix_i\right) = & \sum_{i=0}^n t_i \FldF{f}{v_m}(v_i) = \sum_{k=0}^\ell r_k \FldF{f}{v_m}(x_a) = \ConM{f}{K}\mid_{\Delta^\ell}\left(\sum_{k=0}^\ell r_k x_k\right)\\
= & \sum_{j=0}^m s_j \FldF{f}{v_m}(w_j) = \ConM{f}{K}\mid_{\Delta^\ell\subset\Delta^m}\left(\sum_{j=0}^m s_jw_j\right).
\end{align*}
Obtaining that $\ConM{f}{K}$ is well-defined.
\end{proof}

\begin{lema}
\label{lem:ContinuousConvTrans}
Consider $X$, $Y$, $\mathcal{U}$ and $f$ as in \autoref{def:ConvexTransformation}. $\ConM{f}{}:Y\rightarrow \GRCC{X}$ is continuous, and then $\ConM{f}{}:Y\rightarrow \iGR{X}$ is also continuous. Furthermore, $\ConM{f}{}\simeq f$ in $\iGR{X}$.
\end{lema}

\begin{proof}
By construction, it follows that $\ConM{f}{}$ is continuous in the geometric realization of the simplexes because that map is just the geometric realization from an $n$-simplex. Now, since we have a locally finite closed cover, note that $\ConM{f}{}$ is continuous by 17 A.18 in \cite{Cech_1966}.

The next step is to build the homotopy. Let $H:S^n\times I\rightarrow \iGR{X}$ be the map
\begin{align*}
H(y,t)\coloneqq \left\lbrace\begin{array}{ll}
\ConM{f}{}(y) & \text{ if }t\in [0,1) \\
f(y)  & \text{ if }t=1.
\end{array}\right.
\end{align*}
We claim that $H$ is continuous at any $(y,t)$ with $t<1$.
To prove it, we take a radius small enough such that $\Ball{(y,t)}{r}{S^n\times I}$ does not meet the set $\cap S^n\times \{1\}$, and then
\[H(\Ball{(y,t)}{r}{S^n\times I}) = \Ball{(y,t)}{r}{S^n\times I}.\]
Thus, by continuity in $\ConM{f}{}$, we obtain the continuity at any point of $Y\times [0,1)$.

On the other hand, take $(y,1)$.
By the hypothesis in the open cover and $f$, the image of every simplex $\Delta$ in $Y$ is inside an open set $U$ such that for all $y,y'\in U$ satisfies that $f(y)\in c_{\iGR{X}}(f(y'))$.
Thus, every pair of points $z,z'\in f(\Delta)$ are neighbors as well as they are neighbors the elements of the geometric realization of that simplex (a.k.a $\ConM{f}{}$ of them).
Therefore, since $f(y)\in X\subset \iGR{X}$, for every simplex $\{\Delta_1,\ldots,\Delta_\ell\}$ such that $y\in \Delta_i$ we have that $\cup_{i=1}^\ell \ConM{f}{}(|\Delta_i|)\subset c_{\iGR{X}}(f(y))$.
Furthermore, choosing $r>0$ such that $f(\Ball{y}{r}{Y}) \subset c_{\iGR{X}}(f(y))$ (because $f$ is continuous at $y$), we have that
\begin{align*}
	\ResF{H}{\Delta_i\times I}(\Ball{(y,1)}{r}{S^n\times I}\cap (\Delta_i\times I)) \subset &  \ConM{f}{}(|\Delta_i|)\cup f(\Ball{y}{r}{Y}) \subset \\
	& c_{\iGR{X}}(f(y))
\end{align*}
Thus, $\ResF{H}{\Delta_i\times I}$ is continuous at $(y,1)$ for every $i\in\{1,\ldots,\ell\}$. Hence $H$ is continuous at $(y,1)$. We conclude that $H$ is continuous, and then $\ConM{f}{}\simeq f$ in $\iGR{X}$.
\end{proof}

\begin{observacion}
	Although the map $\ConM{f}{}$ is built as a continuous map from $Y$ to $\GRCC{X}$, we will use the notation $\ConM{f}{}$ both as a map with codomain $\iGR{X}$ and as a map with codomain $\GRCC{X}$ to reduce the notation of composing the identity.
\end{observacion}

\section{Weak Homotopy Equivalence of $(X,\lambda)$ and $\GRCC{X,\lambda}$}
\label{sec.WHEXVRX}

	Remember that the inclusion map $\Id{\uGRCC{X}}: \uGRCC{X}\to \uGRCC{X}$ is continuous from $\GRCC{X}$ to $\iGR{X}$ because $\GRCC{X}$ has a finer structure than the one of $\iGR{X}$.
	From this point on, we will denote the inclusion map between these spaces as $\Id{\uGRCC{X}}$.
	In this section, we will prove that the homomorphism $\pi_n(\Id{\mathbf{KR}}):\pi_n(\GRCC{X})\rightarrow \pi_n(\iGR{X})$ is an isomorphism of groups.
	
	The core of the proof is showing that the functions $S^n \to \iGR{X}$ and $S^n \times I \to \iGR{X}$ are homotopic to maps whose images are in $X\subset \iGR{X}$. These maps have the property that for every point $x$ in the domain, there exists a neighborhood $U$ of $x$ such that every pair of points in $U$ are neighbors.
	Then, we can prove that there is a triangulation of the domain such that the latter map is homotopic to a realization of a simplicial map from this triangulation to $\iGR{X}$.

We denote the set $\{ x\in\mathbb{R}^{n+1} \mid |x|=1 \}$ by $S^n$ considering the topological structure induced by the metric of the geodesic distance metric $d: S^n \times S^n \to \R$  using the standard Riemannian metric. 

Previous to prove that $\Id{\mathbf{KR}}$ is an isomorphism, we first recall the following classical result about compact metric spaces in general.

\begin{lema}[Lebesgue's number lemma; \cite{Willard_2004}, 22.5]
\label{lem:LebesgueNumber}
Let $Y$ be a compact metric space and $\mathcal{U}$ an open cover of $Y$. Then, there exists a number $\delta>0$ such that every subset of $Y$ with a diameter less than $\delta$ is contained in some member of the cover.
\end{lema}

Remember that the numbers that satisfy the condition in \autoref{lem:LebesgueNumber} are called \emph{Lebesgue numbers}. We are now ready to prove the first part of our main result: The group homomorphism $\pi_n(\Id{\textbf{KR}})$ is an epimorphism.

\begin{teorema}
\label{prop:PiNIsEpimorphism}
Let $X$ be a graph. Every map $f:S^n\rightarrow \iGR{X}$ is homotopic to the composition of the identity $\Id{\mathbf{KR}}:\GRCC{X}\rightarrow \iGR{X}$ and some map $g:S^n\rightarrow \GRCC{X}$. This implies that $\pi_n(\Id{\mathbf{KR}})$ is an epimorphism.
\end{teorema}

\begin{proof}
Let $f:(S^n,y_0)\rightarrow (\iGR{X},v_\circledast)$ a continuous based map with basepoints $y_0 \in S^n, v_\circledast\in X$. Since $S^n$ is compact and connected, \autoref{lem:DiscretizationOfAMapConn} implies that $\DMod{X}[f]\simeq f$ in $\iGR{X}$ and $Im(\DMod{X}[f])$ is a finite number of vertices, say $\{v_0,\ldots,v_m\}$. 

Since $\DMod{X}[f]$ is continuous and the set $\{y_0\}$ is compact, we have the flooding discrete modification of $f$ until $i$, $\FldF{f}{v_i}$, of \autoref{def:FldDiscMod}.
It is
\[\FldF{f}{v_i} = \Fld{\cdots \Fld{\Fld{\DMod{X}[f]}{B_0}{v_0}}{B_1}{v_1} \cdots}{B_m}{v_i},\]
where
\[B_i \coloneqq \bigcup_{} \{\Ball{x}{r_{x,i}}{S^n} \mid x \in (\FldF{f}{v_{i-1}})^{-1}(v_i) \}.\]
for some $r_{x,i}$. 
In \autoref{obs:BuildFldDisMod}, we note that
\[f\simeq \FldF{f}{v_m} \text{ and } \{v_\circledast\} = \FldF{f}{v_m}(y_0).\]

Thus, by \autoref{lem:SimplexIsInsideANeighborhood}, for every $y\in S^n$, there exists $r_y>0$ such that for every $w,w'\in \FldF{f}{v_m}(\Ball{y}{r_y}{S^n})$ we have that $w\in c_{\iGR{X}}(w')$. We denote this open cover of $S^n$ as $\mathcal{U}$.

We take a classical triangulation of $S^n$ with $y_0$ as the vertex of some simplices.
Now we apply the barycentric subdivision until we have that every simplex has diameter less than the Lebesgue's number of the cover $\mathcal{U}$, and call that triangulation $K$.
Then, by \autoref{lem:LebesgueNumber} and \autoref{lem:SimplexIsInsideANeighborhood}, the triangulation satisfies the hypothesis of \autoref{lem:ContinuousConvTrans}, and therefore $\ConM{\FldF{f}{v_m}}{K}\simeq \FldF{f}{v_m}$ in $\iGR{X}$. For its definition, $\ConM{\FldF{f}{v_m}}{K}$ is also a continuous map from $S^n$ to $\GRCC{X}$, having the following commutative diagram
\begin{align*}
	\xymatrix{
		& & &(\GRCC{X},v_\circledast) \ar[d]^{\Id{\mathbf{KR}}} \\
		(S^n,y_0) \ar@{->}[rrr]_{\ConM{\FldF{f}{v_m}}{K}} \ar[rrru]^{{\ConM{\FldF{f}{v_m}}{K}}} & & & (\iGR{X},v_\circledast) 
	}
\end{align*}

In conclusion, $f\simeq \FldF{f}{v_m}\simeq \ConM{\FldF{f}{v_m}}{K}$ in $\iGR{X}$.

To prove that this implies that $\pi_n(\Id{\mathbf{KR}})$ is an isomorphism, consider a map $f:(S^n,y_0)\to (\iGR{X},x_0)$. Then $\ConM{\FldF{f}{v_m}}{K}\simeq f \operatorname{rel}\{y_0\}$, $\ConM{\FldF{f}{v_m}}{K}:(Y,y_0)\to (\GRCC{X},x_0)$ is continuous and
\[\pi_n(\Id{\mathbf{KR}})([\ConM{\FldF{f}{v_m}}{K}]) = [\Id{\mathbf{KR}} \circ \ConM{\FldF{f}{v_m}}{K}] = [\ConM{\FldF{f}{v_m}}{K}] = [f].\]
Concluding that $\pi_n$ is group epimorphism.
\end{proof}

To prove $\pi_n(\Id{\mathbf{KR}})$ is a monomorphism, we first need two additional technical results. Following the notation of \cite{Spanier_1966}, we denote the barycentric subdivision of a simplicial complex by $\operatorname{sd}$, and we write $\operatorname{sd}^n$ when barycentric subdivision is applied $n$ times.

\begin{definicion}
\label{def:SimplicialApprox}
Let $K_1$ and $K_2$ be simplicial complexes and let $f:|K_1|\rightarrow |K_2|$ be a continuous map. A simplicial map $\varphi:K_1\rightarrow K_2$ is called a \emph{simplicial approximation} to $f$ if $f(\alpha)\in |s|$ implies that $|\varphi|(\alpha)\in |s|$ for $\alpha\in |K_1|$ and $s_2\in K_2$.
\end{definicion}

\begin{teorema}[Simplicial-approximation Theorem; \cite{Spanier_1966}, 3.4.8]
\label{theo:SimplicialApproxTheo}
Let $(K_1,L_1)$ be a finite simplicial pair and let $f:(|K_1|,|L_1|)\rightarrow (|K_2|,|L_2|)$ be a map. There exists an integer $N$ such that if $n\geq N$ there is a simplicial approximation $(\operatorname{sd}^n K_1,\operatorname{sd}^n L_1)\rightarrow (K_2,L_2)$ to $f$.
\end{teorema}

Thus, for examples, considering the canonical triangulation $K_n$ of $S^n$, there is a non-negative integer $m$ and a simplicial map $\operatorname{sd}^m K_n \to \CC{X}$ such that the geometric realization of the map represents $\beta \in \pi_k(\GRCC{X},v)$.

The following proposition shows that the convex modification is homotopic to the convex modification of the barycentric subdivision in $c_{\GRCC{X}}$.
This result is essential in proving that $\pi_n(\Id{\textbf{KR}})$ is a monomorphism. As we will see, it also imposes some technical conditions in the selection of balls when applying the floods in the proof of \autoref{lem:FloodPreservingACompact}.

\begin{proposicion}
\label{prop:BaryTransOfTheBaryPart}
Let $Y$ be a triangularizable compact metric space with finite triangulation $K$ and let $f: Y=|K|\rightarrow X \subset \iGR{X}$ be such that the image of every simplex in the triangulation $K$ is a simplex in $\iGR{X}$. Then $\ConM{f}{K}$ and $\ConM{f}{\operatorname{sd}K}:Y=|\operatorname{sd}K|\rightarrow \GRCC{X}$ are homotopic.
\end{proposicion}

\begin{proof}
Since $K$ satisfies the condition over the open cover $\mathcal{U}$ to define $\ConM{f}{K}$ in \autoref{def:ConvexTransformation}, then $\operatorname{sd}K$ also satisfies the condition because for every $\sigma'\in\operatorname{sd}K$ there exists a $\sigma\in K$ such that $|\sigma'|\subset |\sigma|$.
Thus $\ConM{f}{\operatorname{sd}K}$ is well defined.

Consider a $n$-simplex $\Delta = \{x_0,\ldots,x_n\}$ in the triangulation $\operatorname{sd}K$ such that $x_0<x_1<\ldots<x_n$ and describe the interval $I$ as the geometric realization of the simplex $\{0,1\}$.
Recall that $|\Delta|\times I$ is a triangularizable metric space with $(n+1)$-simplexes such that if $(w_i,w'_i)\in \Delta\times \{0,1\}$, then either $w_i<w_{i+1}$ and $w'_i=w'_{i+1}$ or $w_i=w_{i+1}$ and $w'_i<w'_{i+1}$.
Let's call this triangulation $K_I$.

Consider the restrictions $\ResF{\ConM{f}{K}}{|\Delta|}: |\Delta|\rightarrow \GRCC{X}$ and $\ResF{\ConM{f}{\operatorname{sd}K}}{|\Delta|}: |\Delta|\rightarrow \GRCC{X}$.
By definition
\begin{align*}
\ConM{f}{\operatorname{sd}K} \mid_{|\Delta|}\left( \sum_{i=0}^n t_ix_i \right) \coloneqq \sum_{i=0}^n t_if(x_i).
\end{align*}
Now, for every $(n+1)$-simplex $\{(w_0,w'_0),\dots, (w_{n+1},w'_{n+1})\}$ the triangulation in $K_I$ of $|\Delta|\times I$, define $H_\Delta:|\Delta|\times I\rightarrow \GRCC{X}$ as
\begin{align*}
H_\Delta\left( \sum_{i=0}^{n+1} t_i (w_i,w'_i) \right) \coloneqq \sum_{i=0}^{n+1} t_i H_\Delta(w_i,w'_i)
\end{align*}
where $H_\Delta(w_i,0)=\ConM{f}{K}(w_i)$ and $H_\Delta(w_i,1)= \ConM{f}{\operatorname{sd }K}(w_i)$.
The map $H_\Delta$ is continuous from $|\Delta|\times I$ to $\GRCC{X}$ because it is continuous in every $(n+1)$-simplex, being a linear map when we restrict it to those $(n+1)$-simplices.

Define $H$ such that $H\mid_{|\Delta|\times I}= H_\Delta$ for every $\Delta$ in $\operatorname{sd}K$.
The map $H$ is continuous because is continuous in $|\Delta|\times I$ and by 17 A.18 in \cite{Cech_1966}. The same construction concludes that the map  $\ConM{f}{\operatorname{sd}K}:|\operatorname{sd}K|\to \GRCC{X}$ is continuous.
By the definition of $H$, we obtain that $H(y,0)=\ConM{f}{K}(y)$ and $H(y,1)=\ConM{f}{\operatorname{sd }K}(y)$.
Hence $\ConM{f}{K}(y)\simeq \ConM{f}{\operatorname{sd }K}$ in $\GRCC{X}$.
\end{proof}

Now we are ready to prove that $\pi_n(\Id{\mathbf{KR}})$ is a monomorphism of groups. 

\begin{proposicion}
\label{prop:PiNIsAMono}
Let $X$ be a graph. $\pi_n(\Id{\mathbf{KR}}):\pi_n(\GRCC{X})\rightarrow \pi_n(\iGR{X})$ is a monomorphism.
\end{proposicion}

\begin{proof}
Let $v_\circledast$ be the base point of $\GRCC{X}$ and denote by $\hat{v}_{\circledast}$ the constant map from $(S^n,y_0)$ to $(\iGR{X},v_\circledast)$ such that $\hat{v}_\circledast(y)=v_\circledast$ for every $y\in S^n$.

Since $\pi_n(\Id{\mathbf{KR}})$ is a homomorphism of groups, it is enough to prove that its kernel has only the identity element.
Suppose that $f:(S^n,y_0)\rightarrow (\GRCC{X},v_\circledast)$ be such that $\Id{\mathbf{KR}}\circ f:(S^n,y_0)\rightarrow (\iGR{X},v_\circledast)$ is homotopic to $\hat{v}_{\circledast}$ in $\iGR{X}$.
Then there is a map $H:S^n\times I\to \iGR{X}$ such that $H(\{y_0\}\times I)=v_\circledast$, $H(y,0)=\Id{\mathbf{KR}} \circ f(y)$ and $H(y,1)=\hat{v}_\circledast(y)=v_\circledast$.
Our goal is to build a homotopy $\hat{H}:S^n \times I \to \GRCC{X}$ from $f$ to $\hat{v_\circledast}$.

Without loss of generality, we can consider $f$ to be a simplicial-approximation to itself by \autoref{theo:SimplicialApproxTheo}. Let $K$ be the triangulation on $S^n$ used in this simplicial approximation. Since we only need a homotopy from $f$ to $\hat{v}_\circledast$, it is enough to consider a map which we know to be homotopic (in $\iGR{X}$) to the original $f$. For this reason, we consider $f$ as $\ConM{\FldF{f}{v_m}}{K}$.

Since $S^n\times I$ is compact, \autoref{lem:DiscretizationOfAMapConn} implies that $\DMod{X}[H]:S^n\times I\rightarrow X\subset \iGR{X}$ satisfies that $\DMod{X} H\simeq H$ and $Im(\DMod{X}[H])$ is finite, say $\{v_0,\ldots,v_m\}$.

Giving that $\DMod{X}[H]$ is continuous in $\iGR{X}$ and $\{y_0\}\times I$ is compact, the flooding discrete modification of $H$ until $v_i$ as in \autoref{def:DiscModMap} is defined. It is
\[\FldF{H}{v_i} = \Fld{\cdots \Fld{\Fld{\DMod{X}[H]}{B_0}{v_0}}{B_1}{v_1} \cdots}{B_m}{v_i},\]
where
\[B_i \coloneqq \bigcup_{} \{\Ball{y}{r_{(y,t),i}}{S^n} \mid (y,t) \in (\FldF{H}{v_{i-1}})^{-1}(v_i) \}.\]
for some $r_{(y,t),i}$ such that $\Ball{(y,t)}{r_{(y,t),i}}{S^n\times I}$ does not meet $(\{y_0\}\times I)$ if $v_i\neq v_\circledast$. Moreover, we can ask for $r_{(y,t),i}$ such that
\begin{itemize}[wide,labelwidth=!, labelindent=0pt] 
	\item if $t=1$, $\Ball{(y,t)}{r_{(y,t),i}}{S^n\times I}\cap S^n\times \{0\}=\varnothing$.
	\item if $0<t<1$, $\Ball{(y,t)}{r_{(y,t)}}{S^n\times I}\cap S^n\times \{0,1\}=\varnothing$.
	\item if $t=0$, $\Ball{(y,t)}{r_{(y,t)}}{S^n\times I}\cap S^n\times \{1\}=\varnothing$ and $\Ball{(y,t)}{r_{(y,t)}}{S^n\times I}$ does not meet any vertex of the triangulation of $S^n$ except for $(y,t)$ itself (in the case it is a vertex). This last condition is possible because we have just a finite number of simplexes.
\end{itemize}

By the choose of the $r_{(y,t),i}$, $H\simeq \FldF{H}{v_m}$ and $\FldF{H}{v_m} = \{v_\circledast\}$. \autoref{lem:SimplexIsInsideANeighborhood} implies that for every $(y,t)\in S^n\times I$ there exists $\rho_{(y,t)}>0$ such that for every $w,w'\in \FldF{H}{v_m}(\Ball{(y,t)}{\rho_{(y,t)}}{S^n\times I}$ we have that $w\in c_{\iGR{X}}(w')$. Let's denote by
\[\mathcal{U} = \{ \Ball{(y,t)}{\rho_{(y,t)}}{S^n\times I} \mid (y,t)\in S^n\times I \}\] this open cover of $S^n\times I$.

Considering the triangulation $K$ of $S^n$ and $I$ as the geometric realization of $\{0,1\}$, take the barycentric subdivision of both $K$ and $I$ until the diameter of the simplices in $\operatorname{sd}^k K \times \operatorname{sd}^k \{0,1\}$ is less than the Lebesgue number of $\mathcal{U}$.
Then, the triangulation $\operatorname{sd}^k K \times \operatorname{sd}^k \{0,1\}$ satisfies the hypothesis of \autoref{lem:ContinuousConvTrans} over the open cover $\mathcal{U}$.
Thus $\ConM{\FldF{H}{v_m}}{\operatorname{sd}^k K \times \operatorname{sd}^k \{0,1\}}$ is continuous in $\GRCC{X}$, $\ConM{\FldF{H}{v_m}}{\text{sd}^k K \times \operatorname{sd}^k \{0,1\}}(y,1) = \hat{v}_\circledast(y)=v_\circledast$, and $\ConM{\FldF{H}{v_m}}{\operatorname{sd}^k K \times \operatorname{sd}^k \{0,1\}}(y,0) = \ConM{\FldF{f}{v_m}}{\operatorname{sd}^k K}(y)$.
Hence $\ConM{\FldF{f}{v_m}}{\operatorname{sd}^k K} \simeq \hat{v}_\circledast$ in $\GRCC{X}$.

Finally, since the points of the vertices of the triangulation were not changed, \autoref{prop:BaryTransOfTheBaryPart} implies the following sequence of homotopies in $\GRCC{X}$
\begin{align*}
\ConM{\FldF{f}{v_m}}{\operatorname{sd}^{k} K} \simeq \ConM{\FldF{f}{v_m}}{\operatorname{sd}^{k-1} K} \simeq \ldots \simeq \ConM{\FldF{f}{v_m}}{\operatorname{sd} K} \simeq \ConM{\FldF{f}{v_m}}{K} = f
\end{align*}
Concluding that $f\simeq \hat{v}_\circledast$ in $\GRCC{X}$, as desired.
It is, if $\pi_n(\Id{\textbf{KR}})([f]) = [\hat{v}_\circledast]$, then $[f] = [\hat{v}_\circledast]\in \pi_n(\GRCC{X},v_\circledast)$, proving that $\pi_n(\Id{\textbf{KR}})$ is a monomorphism.
\end{proof}

\begin{teorema}
\label{Theo:XandRealizationOfVRSameWHT}
Let $G = (X, E)$ be a graph, and let $(X,c_G)$ be the induced \v{C}ech closure space. Then $\GRCC{X}$ are weakly homotopy equivalent as $(X,c_X)$.
\end{teorema}

\begin{proof}
	Given the strong homotopy equivalence in \autoref{prop:SHEBetweenGraphandRGraph} between the  maps $\Id{X}:X\rightarrow \iGR{X}$ and $\DMod{X}:\iGR{X}\rightarrow X$, we obtain that $\pi_n(\DMod{X}):\pi_n(\iGR{X})\rightarrow \pi_n(X)$ is an isomorphism.

	In addition, by \autoref{prop:PiNIsEpimorphism} and \autoref{prop:PiNIsAMono}, we have that $\pi_n(\Id{\mathbf{KR}}):\pi_n(\GRCC{X})\rightarrow \pi_n(\iGR{X})$ is also an isomprhism. Thus, we obtain that \[\pi_n(\DMod{X}\circ \Id{\mathbf{KR}}) = \pi_n(\DMod{X})\circ \pi_n(\Id{\mathbf{KR}}):\pi_n(\GRCC{X})\rightarrow \pi_n(X)\] is an isomorphism.
\end{proof}

\section{Discussion}
\label{sec.App}

In this brief section, we connect the main result of this paper \autoref{Theo:XandRealizationOfVRSameWHT} , to the findings of two other papers.

First, Milićević and Scoville prove in their Theorem 24 of \cite{MilicevicScoville2026} that each finite digraph $X$ is weakly homotopy equivalent to the geometric realization of its directed Vietoris-Rips complex.
They prove this by using a map, $f_X$, from the latter to the former that is similar to the discrete modification defined in the geometric realization of the directed Vietoris-Rips complex.
In contrast, we proved the same result for arbitrary graphs using the intermediate geometric realization.

Second, Adamaszek and Adams proved in \cite{Adamaszek_Adams_2017} that
\[\mathbf{VR}_{<}(S^1,r) \simeq S^{2\ell+1} \text{ if }\frac{\ell}{2\ell+1}< r \leq \frac{\ell+1}{2\ell+3}\]
where $\mathbf{VR}_{<}(S^1,r)$ is the Vietoris-Rips simplicial complex, i.e., its simplices are the finite subsets of $S^1$ of diameter less than $r$.
Considering the graph $(S^1,E)$ such that $\{v,w\}\in E$ if $d(v,w)<r$, we note that $\mathbf{VR}_{<}(S^1,r) = \CC{(S^1,E)}$. Furthermore, if we define $c_{<r}:\PowS{S^1}\to \PowS{S^1}$ such that \[c_{<r}(A) = \{x\in S^1\mid d(x,A)<r\}.\]
Then
\[c_{<r}(A) = \bigcup\{ c_{<r}(x) \mid x\in A\}.\]
Thus, the canonical closure operator for the graph $(S^1,E)$, $c_{(S^1,E)}$ is equal to $c_{<r}$.
Hence, we conclude that
\[\mathbf{VR}_{<}(S^1,r) \simeq \GRCC{(S^1,E)} \simeq (S^1,c_{(S^1,E)}) = (S^1,c_{<r})\]
where the second $\simeq$, from left to right, is the result of the main theorem.

\section*{Acknowledgments}
	I would like to thank my thesis adviser, Antonio Rieser, who has helped with several suggestions about the order of the ideas of this paper to improve the delivery and to reduce the number of particular lemmas for specific cases.
	I am deeply grateful to the AIM Workshop on Discrete and Combinatorial Homotopy Theory, its organizers: Helen Barcelo, Antonio Rieser and Volkmar Welker, and certainly with my one-week teammates: Greg Lupton, Oleg Musin, Nicholas Scoville and Chris Staecker;
	I learned a great deal from them about how to decompose a problem, discuss the ideas and try to simplify them in order to improve the exposition.
	I am grateful to Antonio Rieser for telling me about the existence of \cite{WebJMMNicks} after I showed him the first versions of this work.

\printnomenclature

\bibliographystyle{amsalpha}
\bibliography{all}

\end{document}